# New multivariate central limit theorems in linear structural and functional error-in-variables models

**Yuliya V. Martsynyuk**[*]

*School of Mathematics and Statistics, Carleton University,
1125 Colonel By Drive, Ottawa, ON K1S 5B6, Canada.
e-mail:* `ymartsynyuk@yahoo.com`

**Dedicated to Miklós Csörgő on the occasion of his 75th birthday**

**Abstract:** This paper deals simultaneously with linear structural and functional error-in-variables models (SEIVM and FEIVM), revisiting in this context generalized and modified least squares estimators of the slope and intercept, and some methods of moments estimators of unknown variances of the measurement errors. New joint central limit theorems (CLT's) are established for these estimators in the SEIVM and FEIVM under some first time, so far the most general, respective conditions on the explanatory variables, and under the existence of four moments of the measurement errors. Moreover, due to them being in Studentized forms to begin with, the obtained CLT's are a priori nearly, or completely, data-based, and free of unknown parameters of the distribution of the errors and any parameters associated with the explanatory variables. In contrast, in related CLT's in the literature so far, the covariance matrices of the limiting normal distributions are, in general, complicated and depend on various, typically unknown parameters that are hard to estimate. In addition, the very forms of the CLT's in the present paper are universal for the SEIVM and FEIVM. This extends a previously known interplay between a SEIVM and a FEIVM. Moreover, though the particular methods and details of the proofs of the CLT's in the SEIVM and FEIVM that are established in this paper are quite different, a unified general scheme of these proofs is constructed for the two models herewith.

**AMS 2000 subject classifications:** Primary 60F05, 62J99; secondary 60E07.
**Keywords and phrases:** linear structural/functional error-in-variables model, measurement errors, explanatory variables, domain of attraction of the normal law, slowly varying function, identifiability assumptions, generalized/modified least squares estimator, central limit theorem, multivariate Student statistic, positive definite matrix, Cholesky square root of a matrix, symmetric positive definite square root of a matrix, Lindeberg's condition, generalized domain of attraction of the multivariate normal law, spherically symmetric random vector, full random vector.

Received June 2007.

---

[*]Research supported by a Carleton University Faculty of Graduate Studies and Research scholarship, an NSERC Canada Discovery Grant of M. Csörgő at Carleton University, and an NSERC Postdoctoral Fellowship of Yu.V. Martsynyuk at University of Ottawa.

**Contents**

## 1. Introduction

First, we present linear structural and functional error-in-variables models, assumptions that will be used in each of them, and estimators of unknown parameters of interest under study (cf. respective Sections 1.1–1.3). We then discuss the essence of the main results of this paper in Section 1.4.

### 1.1. Linear structural and functional error-in-variables models (SEIVM and FEIVM)

In the linear structural and functional error-in-variables models (EIVM's) of this paper we observe pairs $(y_i, x_i) \in \mathbb{R}^2$ according to

$$y_i = \beta \xi_i + \alpha + \delta_i, \quad x_i = \xi_i + \varepsilon_i, \qquad (1.1)$$

where $\xi_i$ are unknown explanatory/latent variables, the real-valued slope $\beta$ and intercept $\alpha$ are to be estimated, and $\delta_i$ and $\varepsilon_i$ are unknown measurement error terms/variables, $1 \le i \le n$, $n \in \mathbb{N}$. EIVM (1.1) is also known as a measurement error model, or structural/functional relationship, or regression with errors in variables. It is a generalization of the simple linear regression of form $y_i =$ $\beta \xi_i + \alpha + \delta_i$ in that in (1.1) it is assumed that, in addition to the two variables $\eta := \beta \xi + \alpha$ and $\xi$ being linearly related, now not only $\eta$, but also $\xi$, are observed with respective measurement errors $\delta_i$ and $\varepsilon_i$.

The explanatory variables $\xi_i$ are assumed to be independent identically distributed (i.i.d.) random variables (r.v.'s) that are independent of the error terms when we deal with the structural EIVM (SEIVM) (cf. upcoming condition **(S2)**), and deterministic in case of the functional EIVM (FEIVM).

### 1.2. Assumptions in SEIVM and FEIVM

Both in the SEIVM and FEIVM versions of (1.1), the following assumptions are made on the error terms.

**(A)** $\{(\delta, \varepsilon), (\delta_i, \varepsilon_i), i \geq 1\}$ is a sequence of i.i.d. random vectors with mean zero, $E\delta^4 < \infty$ and $E\varepsilon^4 < \infty$, and with a positive definite covariance matrix

$$\Gamma = \text{Cov}(\delta, \varepsilon) = \begin{pmatrix} \lambda\theta & \mu \\ \mu & \theta \end{pmatrix}. \tag{1.2}$$

**(B)** The r.v.'s $\delta$ and $\varepsilon$ are either both continuous, or one of them is continuous, and the other one is discrete.

In the SEIVM, we suppose that the explanatory variables $\{\xi, \xi_i, i \geq 1\}$ obey **(S1)** and **(S2)** as follows:

**(S1)** $\{\xi, \xi_i, i \geq 1\}$ are i.i.d.r.v.'s in the domain of attraction of the normal law (DAN), i.e., there are constants $a_n$ and $b_n$, $b_n > 0$, for which $(\sum_{i=1}^n \xi_i - a_n)b_n^{-1} \xrightarrow{\mathcal{D}} N(0,1)$, as $n \to \infty$.

**(S2)** $\xi$ is independent of $(\delta, \varepsilon)$.

**Remark 1.1.** Further to the definition of DAN in **(S1)**, it is known that $a_n$ can be taken as $nE\xi$ and $b_n = n^{1/2}\ell_\xi(n)$, where $\ell_\xi(n)$ is a slowly varying function at infinity (i.e., $\ell_\xi(az)/\ell_\xi(z) \to 1$, as $z \to \infty$, for any $a > 0$), defined by the distribution of $\xi$. Moreover, $\ell_\xi(n) = \sqrt{\text{Var}\,\xi} > 0$, if $\text{Var}\,\xi < \infty$, and $\ell_\xi(n) \nearrow \infty$, as $n \to \infty$, if $\text{Var}\,\xi = \infty$. Also, $\xi$ has moments of all orders less than 2, and the variance of $\xi$ is positive, but need not be finite.

**Remark 1.2.** One of the several necessary and sufficient conditions for i.i.d.r.v.'s $\{Z, Z_i, i \geq 1\}$ to be in DAN is commonly associated with O'Brien [19] (for more details see also Remark (iii) in Maller [13], p.194), and it reads as follows: $\max_{1 \leq i \leq n} Z_i^2 \big/ \sum_{i=1}^n Z_i^2 \xrightarrow{P} 0$, $n \to \infty$. In Remark (iv) of Maller [13] it is pointed out that this "negligibility" condition has appeared and played an important role in the asymptotic theory of many stochastic models. In addition to the models listed in [13], $Z \in$ DAN has also been used for CLT's in a simple linear regression (cf. Maller [12]) and, as frequently an optimal, or nearly optimal, condition on the explanatory variables, for various marginal CLT's in the SEIVM (1.1) (cf. Martsynyuk [14, 15, 17]).

For the deterministic $\xi_i$ in the FEIVM (1.1), we assume the following **(F1)**–**(F3)** conditions that, in view of Remarks 1.1 and 1.2, are natural companions to **(S1)**.

**(F1)** $\lim_{n\to\infty} n^{-1} \sum_{i=1}^n \xi_i = m$ and $m$ is finite.
**(F2)** Either $0 < \lim_{n\to\infty} n^{-1} \sum_{i=1}^n \xi_i^2 < \infty$, or $\lim_{n\to\infty} n^{-1} \sum_{i=1}^n \xi_i^2 = \infty$.
**(F3)** $\dfrac{\max_{1\le i \le n} \xi_i^2}{\sum_{i=1}^n \xi_i^2} \to 0, \quad n \to \infty.$

For some further comments on the introduced assumptions on the explanatory variables in the SEIVM (1.1) and FEIVM (1.1), we refer to Remarks 2.10 and 2.11 of Section 2.2.

To ensure identifiability of unknown parameters in the SEIVM and FEIVM versions of (1.1), it is common to make use of some additional side conditions in this regard, usually as conditions on the matrix $\Gamma$ of (1.2) in **(A)**. Here we distinguish three most frequently used identifiability assumptions, namely:

**(1)** the positive ratio of the error variances $\lambda = \operatorname{Var}\delta/\operatorname{Var}\varepsilon$ is known, and the covariance of $\delta$ and $\varepsilon$ is zero, i.e., $\operatorname{cov}(\delta, \varepsilon) = \mu = 0$;
**(2)** $\operatorname{Var}\delta = \lambda\theta$ and $\operatorname{cov}(\delta,\varepsilon) = \mu$ are known, while $\operatorname{Var}\varepsilon = \theta$ is unknown;
**(3)** $\operatorname{Var}\varepsilon = \theta$ and $\operatorname{cov}(\delta,\varepsilon) = \mu$ are known, while $\operatorname{Var}\delta = \lambda\theta$ is unknown.

The identifiability assumption **(1)** corresponds to orthogonal regression estimation in (1.1), while **(3)** is likely to be realistic in many applications (cf. Carroll and Ruppert [1], Carroll *et al.* [2], Cheng and Van Ness [4] and Fuller [6] for some further discussions along these lines).

### 1.3. *Estimators for slope $\beta$, intercept $\alpha$ and unknown error variances*

For further use throughout, for real-valued variables $\{u_i, 1 \le i \le n\}$ and $\{v_i, 1 \le i \le n\}$, we put

$$\overline{u} = \frac{1}{n} \sum_{i=1}^n u_i, \quad s_{i,uv} = (u_i - \overline{u})(v_i - \overline{v}) \quad \text{and} \quad S_{uv} = \frac{1}{n} \sum_{i=1}^n s_{i,uv}.$$

Under each of the identifiability assumptions in **(1)**–**(3)**, we are to estimate the slope $\beta$, intercept $\alpha$, and the respective, typically unknown, error variance $\lambda\theta$, or $\theta$, that is denoted throughout by $\gamma$ for convenience in notations.

When **(1)** is assumed, it has been common to estimate $\beta$ and $\alpha$ with generalized least squares estimators (GLSE's) $\widehat{\beta}_{1n}$ and $\widehat{\alpha}_{1n}$ that are simply derived under **(A)** without assuming the finiteness of the fourth error moments (cf., e.g., Section 3.4 in [4]). For estimating $\gamma = \theta$, a method of moments estimator (MME) $\widehat{\gamma}_{1n}$ is a usual one. Provided that $S_{xy} \neq 0$, these estimators are given

by

$$\widehat{\beta}_{1n} = \text{sign}(S_{xy})\sqrt{\left(\frac{\lambda S_{xx} - S_{yy}}{2S_{xy}}\right)^2 + \lambda} - \frac{\lambda S_{xx} - S_{yy}}{2S_{xy}},$$
$$\widehat{\alpha}_{1n} = \overline{y} - \overline{x}\widehat{\beta}_{1n}, \ \widehat{\gamma}_{1n} = \frac{S_{yy} - 2S_{xy}\widehat{\beta}_{1n} + S_{xx}\widehat{\beta}_{1n}^2}{\lambda + \widehat{\beta}_{1n}^2}, \quad (1.3)$$

where $\text{sign}(S_{xy})$ denotes the sign of $S_{xy}$.

When either **(2)** or **(3)** is assumed, modified least squares estimators (MLSE's) for $\beta$ and $\alpha$ (cf. Section 3.5 in [4]), and MME's for the respective unknown error variances $\gamma = \theta$ and $\gamma = \lambda\theta$ are available. The MLSE's and MME under **(2)** are

$$\widehat{\beta}_{2n} = \frac{S_{yy} - \lambda\theta}{S_{xy} - \mu}, \ \widehat{\alpha}_{2n} = \overline{y} - \overline{x}\widehat{\beta}_{2n} \text{ and } \widehat{\gamma}_{2n} = S_{xx} - \frac{S_{xy} - \mu}{\widehat{\beta}_{2n}}, \quad (1.4)$$

assuming that $S_{xy} - \mu \neq 0$ and $S_{yy} - \lambda\theta > 0$, while those under **(3)** are

$$\widehat{\beta}_{3n} = \frac{S_{xy} - \mu}{S_{xx} - \theta}, \ \widehat{\alpha}_{3n} = \overline{y} - \overline{x}\widehat{\beta}_{3n} \text{ and } \widehat{\gamma}_{3n} = S_{yy} - (S_{xy} - \mu)\widehat{\beta}_{3n}, \quad (1.5)$$

provided that $S_{xx} - \theta > 0$.

### *1.4. Introduction to main results*

In this paper, we revisit the triples of the estimators in (1.3)–(1.5) in the SEIVM and FEIVM versions of (1.1), and prove two joint central limit theorems (CLT's) for each of the triples (cf. Theorem 2.1), under the respectively introduced assumptions on the explanatory variables that are, to the best of our knowledge, the most general ever used so far in this context (cf. Remarks 2.2, 2.3). As to the conditions **(A)** and **(B)** on the error terms here, they seem to be the least restrictive that have been considered in the literature thus far (cf. Remarks 2.2, 2.3, 2.8, 2.9).

Further to the special features of our CLT's in Theorem 2.1, these CLT's are in Studentized forms to begin with and, as a result, are automatically nearly, or completely, data-based. Namely, as compared to the related CLT's for $(\beta, \alpha, \gamma)$ in the literature, the ones in Theorem 2.1 are a priori free of any unknown parameters associated with the explanatory and error variables (cf. Remarks 2.5 and 2.6).

The CLT's of Theorem 2.1 also extend a previously known interplay between a SEIVM and a FEIVM as in Gleser [8]. This extension is due to a synchronized choice of the respective conditions on the explanatory variables in the SEIVM and FEIVM of the present paper (cf. Remark 2.11). Consequently, the CLT's of Theorem 2.1 are universal in form for the SEIVM (1.1) and FEIVM (1.1) (cf. Remark 2.12).

The idea of establishing the joint CLT's of Theorem 2.1 for the estimators in (1.3)–(1.5) has originated from a wish to extend and unify the Studentized marginal CLT's for each of these estimators that are proved, among other things, for the SEIVM (1.1) in Martsynyuk [14, 15, 17], and for the FEIVM (1.1) in Martsynyuk [15, 16], under nearly the same respective model assumptions as those in Theorem 2.1. When establishing the multivariate CLT's of the present paper, it was important for us to preserve and build on the new assumptions on the explanatory variables that had first been introduced and used in the SEIVM (1.1) and FEIVM (1.1) in [14, 15] (cf. Remarks 2.10 and 2.11 on the crucial roles of **(S1)** and **(F1)**–**(F3)** in [14, 15]). However, it was even more desirable that the CLT's for $(\beta, \alpha, \gamma)$ should be of suitable Studentized forms, like their marginal predecessors, and so that they would also be universal in form both for the SEIVM (1.1) and FEIVM (1.1).

Theorem 2.1 would have remained a wishfull thinking only if not for the auxiliary results in Section 3 that bridge the context of the SEIVM (1.1) and FEIVM (1.1) with recent advances in Studentization of random vectors by a matrix, the generalized domain of attraction of the multivariate normal law (GDAN), and the domain of attraction of the univariate normal law (DAN). For the sake of providing convenient reference in Section 3, we summarize some of these advances in the subsidiary Section 4. Among the auxiliary results of Section 3, the key Theorems 3.1 and 3.2 of Sections 3.1 and 3.2, respectively, namely, a CLT for a multivariate Student statistic that is based on independent but not necessarily identically distributed random vectors that satisfy the Lindeberg condition, and a special characterization of GDAN, may also be of interest beyond the scope of the present paper. Also, the auxiliary CLT's of Section 3.3 are rather versatile, and can also be used to prove multivariate CLT's for estimators other than (1.3)–(1.5) in the SEIVM (1.1) and FEIVM (1.1), and in the respective no-intercept versions of these models, where $\alpha$ is assumed to be zero (cf. Remarks 3.4 and 3.5). Although the particular methods and details of the proofs of Theorem 2.1 for the SEIVM and FEIVM versions of (1.1) are fundamentally different, a unified general scheme of these proofs is constructed for the two models (cf. Sections 3.3 and 3.4).

This paper is based on parts of the author's Ph.D. thesis Martsynyuk [15], written under the supervision of Miklós Csörgő, and on parts of Martsynyuk [14].

## 2. Main results

### 2.1. Main results with remarks

The (a) and (b) parts of Theorem 2.1, namely the Studentized CLT's for the triples of the estimators in (1.3)–(1.5) in the SEIVM (1.1) and FEIVM (1.1), constitute the main results of this paper.

In the sequel, all vectors are row-vectors. For vectors $Z_1, \cdots, Z_n$ and $W_1, \cdots, W_n$ in $\mathbb{R}^d$, $d \geq 1$, introduce vector $\overline{Z} = n^{-1} \sum_{i=1}^n Z_i$ and matrix $V_{ZW} = (n-1)^{-1} \sum_{i=1}^n s_{i,ZW}$, where $s_{i,ZW} = (Z_i - \overline{Z})^T (W_i - \overline{W})$ and $(Z_i - \overline{Z})^T$ is the transpose of $Z_i - \overline{Z}$. For a positive definite matrix $A$ ($A > 0$), notation $A^{1/2}$ stands both for the (left) Cholesky and symmetric positive definite square roots of $A$. We recall that the (left) Cholesky square root $A^{1/2}$ of $A > 0$ is the uniquely existing lower triangular matrix with positive diagonal elements that is such that $A^{1/2}(A^{1/2})^T = A$. Clearly, it is invertible. As to the symmetric positive definite square root $A^{1/2}$ of matrix $A > 0$, the latter exists and satisfies $(A^{1/2})^2 = A$, where $A^{1/2} = (A^{1/2})^T$. By definition, $A^{T/2} = (A^{1/2})^T$, $A^{-1/2} = (A^{1/2})^{-1}$ and $A^{-T/2} = (A^{-1/2})^T$. Notation $\mathrm{diag}(\cdot, \cdots, \cdot)$ stands for a block-diagonal matrix, where in the brackets square matrix blocks that are on its diagonal are listed.

**Theorem 2.1.** *Let assumptions*

$$\begin{cases} \textbf{(A)}, \textbf{(B)}, \textbf{(S1)} \text{ and } \textbf{(S2)}, \text{ in the SEIVM } (1.1), \\ \textbf{(A)}, \textbf{(B)} \text{ and } \textbf{(F1)}-\textbf{(F3)}, \text{ in the FEIVM } (1.1), \end{cases}$$

*be satisfied. Suppose also that the identifiability condition in* **(1)**–**(3)** *that is appropriate for the estimators in hand is valid. Let*

$$U(j,n) = \begin{cases} 2S_{xy} & , \text{ if } j = 1, \\ S_{xy} - \mu & , \text{ if } j = 2, \\ S_{xx} - \theta & , \text{ if } j = 3, \end{cases} \quad L(j,n) = \begin{cases} n^{-1}(n-2)(\lambda + \widehat{\beta}_{1n}^2) & , \text{ if } j = 1, \\ 1 & , \text{ if } j = 2, \\ 1 & , \text{ if } j = 3, \end{cases} \quad (2.1)$$

$$u_i(j,n,\beta) = \begin{cases} \dfrac{-2\beta^2}{\lambda + \beta^2} \left( \lambda s_{i,xx} - s_{i,yy} - \dfrac{\lambda - \beta^2}{\beta} s_{i,xy} \right) & , \text{ if } j = 1, \\ (s_{i,yy} - \lambda\theta) - \beta(s_{i,xy} - \mu) & , \text{ if } j = 2, \\ (s_{i,xy} - \mu) - \beta(s_{i,xx} - \theta) & , \text{ if } j = 3, \end{cases} \quad (2.2)$$

$$v_i(j,n,\beta) = y_i - \alpha - \beta x_i - \dfrac{\overline{x}}{U(j,n)} u_i(j,n,\beta), \quad 1 \leq j \leq 3, \quad (2.3)$$

$$w_i(j,n,\beta) = \begin{cases} (s_{i,yy} - \lambda\theta) - 2\beta(s_{i,xy} - \mu) + \beta^2(s_{i,xx} - \theta) & , \text{ if } j = 1, \\ \beta^{-2} w_i(1,n,\beta) & , \text{ if } j = 2, \\ w_i(1,n,\beta) & , \text{ if } j = 3, \end{cases} \quad (2.4)$$

$$z_i(j,n,\beta) = \Big(u_i(j,n,\beta), v_i(j,n,\beta), w_i(j,n,\beta)\Big), \quad 1 \leq j \leq 3. \quad (2.5)$$

*Then, for $1 \leq j \leq 3$, as $n \to \infty$, the following* CLT's *hold true:*

(a) $\sqrt{n} \Big( U(j,n)(\widehat{\beta}_{jn} - \beta), \widehat{\alpha}_{jn} - \alpha, L(j,n)(\widehat{\gamma}_{jn} - \gamma) \Big) V_{z(j,n,\beta)z(j,n,\beta)}^{-T/2} \xrightarrow{\mathcal{D}} N(0, I_3)$;

(b) $\sqrt{n} \Big( U(j,n)(\widehat{\beta}_{jn} - \beta), \widehat{\alpha}_{jn} - \alpha, L(j,n)(\widehat{\gamma}_{jn} - \gamma) \Big) V_{z(j,n,\widehat{\beta}_{jn})z(j,n,\widehat{\beta}_{jn})}^{-T/2} \xrightarrow{\mathcal{D}} N(0, I_3)$.

**Remark 2.1.** Matrices $V_{z(j,n,\beta)z(j,n,\beta)}^{-T/2}$ and $V_{z(j,n,\widehat{\beta}_{jn})z(j,n,\widehat{\beta}_{jn})}^{-T/2}$ in the respective (a) and (b) parts of Theorem 2.1 are well defined both in the SEIVM and FEIVM. Indeed, for $j = 3$, on combining Remark 4.1, (3.80) and (3.87) with $p_i(n)$ of (3.86), we conclude that $V_{z(3,n,\beta)z(3,n,\beta)} > 0$ and $V_{z(3,n,\widehat{\beta}_{jn})z(3,n,\widehat{\beta}_{jn})} > 0$ on sets whose probabilities converge to one, as $n \to \infty$. We note that matrices $B_n$ in (3.80) are invertible in view of (3.72) and Remark 3.3. Similar arguments apply in the case of $j = 1$ or $j = 2$.

**Remark 2.2.** In view of Remark 1.1, the assumption $\xi \in$ DAN of **(S1)** in the SEIVM (1.1) is weaker than the following one:

$$E\xi = m \text{ and } \text{Var}\,\xi = E\xi^2 - (E\xi)^2 =: M - m^2 > 0, \text{ with } M < \infty. \quad (2.6)$$

Also, conditions **(F1)**–**(F3)** on the explanatory variables in the FEIVM (1.1) are less restrictive than those of

$$\overline{\xi} \to m \text{ and } \overline{\xi^2} \to M, \text{ as } n \to \infty, \text{ and } M - m^2 > 0, \text{ with finite } m \text{ and } M. \quad (2.7)$$

Indeed, convergence of $\overline{\xi^2}$ to the finite positive limit $M$ implies that $n^{-1}\xi_n^2 = \overline{\xi^2} - n^{-1}(n-1)(n-1)^{-1}\sum_{i=1}^{n-1}\xi_i^2 \to 0$, $n \to \infty$, and the latter convergence, via a proof by contradiction, leads to $n^{-1}\max_{1\le i\le n}\xi_i^2 \to 0$, $n \to \infty$, and thus, to **(F3)**. While conditions (2.6) and (2.7) have commonly been used for the CLT's in the SEIVM's and FEIVM's in the literature so far, our **(S1)** and the group of assumptions of **(F1)**–**(F3)** are believed to be new respective assumptions for these models. More precisely, **(S1)** and **(F1)**–**(F3)** have first been introduced respectively for the SEIVM (1.1) and FEIVM (1.1) in Martsynyuk [14, 15] (cf. also Remarks 2.10 and 2.11), and have not yet been used in EIVM's by other authors. As to the conditions on the error terms in **(A)** and **(B)** here, they seem to be the least restrictive that have been considered for CLT's thus far (we will elaborate further on the assumption **(B)** in Remarks 2.8, 2.9).

**Remark 2.3.** In the literature, the vectors of the estimators $(\widehat{\beta}_{jn}, \widehat{\alpha}_{jn}, \widehat{\gamma}_{jn})$, $1 \le j \le 3$, are known to be $\sqrt{n}$-asymptotically normal (cf., e.g., Gleser [7, 8] for $j = 1$, Cheng and Van Ness [3] for $j = 2$ and 3), under the respective identifiability assumption in **(1)**–**(3)**, **(A)**, and (2.6) and **(S2)** in the SEIVM, or (2.7) in the FEIVM, and, in case of [3], also under the condition that $\xi$, $\delta$ and $\varepsilon$ are indepedent normal r.v.'s in the SEIVM, and that $\delta$ and $\varepsilon$ are independent normal r.v.'s in the FEIVM. In view of Remark 2.2, if (2.6) and (2.7) are not assumed, then the CLT's of Theorem 2.1 appear to be first time ones. Otherwise, under (2.6) (in the SEIVM), or (2.7) (in the FEIVM), subject to the above mentioned additional normality conditions in [3] and comments on **(B)** in upcoming Remarks 2.8 and 2.9, the CLT's of [3, 7, 8] follow respectively from those in Theorem 2.1. For the proof of this statement, we refer to the very end of Section 3.4.

**Remark 2.4.** Theorem 2.1 is a unifying generalization of the author's Studentized marginal CLT's for $\widehat{\beta}_{jn}$, $\widehat{\alpha}_{jn}$ and $\widehat{\gamma}_{jn}$, $1 \le j \le 3$, that are established for the SEIVM (1.1) and FEIVM (1.1) respectively in [14, 15, 17] and [15, 16]. Also, Theorem 2.1 is proved under nearly the same assumptions as those used for the marginal CLT's in [14, 15, 16, 17].

**Remark 2.5.** Due to them being in Studentized forms to begin with, the CLT's of Theorem 2.1 are a priori free of any unknown parameters (such as moments) of the distribution of $(\delta, \varepsilon)$ depend only on the error moments that are assumed to be known according to the corresponding identifiability assumption in **(1)**–**(3)**), and do not contain any parameters associated with $\xi$ in the SEIVM and $\{\xi_i, i \geq 1\}$ in the FEIVM. In addition, the CLT's of the (b) part of Theorem 2.1 are completely data-based. Hence, the latter CLT's are readily applicable to constructing large-sample approximate confidence regions for $(\beta, \alpha, \gamma)$.

**Remark 2.6.** A priori Studentized forms, and the corresponding features described in Remark 2.5, of our CLT's make them new even under the stronger conditions in (2.6) and (2.7) that were used in the CLT's of [3, 7, 8]. Indeed, as opposed to Theorem 2.1, the expression (the same for the SEIVM and FEIVM) for the covariance matrix of the asymptotic normal distribution of $(\widehat{\beta}_{1n}, \widehat{\alpha}_{1n}, \widehat{\gamma}_{1n})$ that is due to Gleser [7, 8] is complicated, and in addition to the unknown parameters $\beta$, $m$ and $M$, where $m$ and $M$ are as in (2.6) and (2.7), it involves typically unknown cross-moments and moments of order $\leq 4$ of the error terms that are hard-to-estimate from data. Then, in order to be able to estimate the covariance matrix of the CLT in [7] (in the FEIVM), it is additionally assumed that the moments of $\delta$ and $\varepsilon$ are like those of two independent normal r.v.'s. Consequently, the covariance matrix of the latter CLT becomes much simpler in form and contains only the unknown, but estimable $\beta$, $m$, $M$ and $\gamma = \lambda\theta$. As to the respective asymptotic covariance matrices of $(\widehat{\beta}_{jn}, \widehat{\alpha}_{jn}, \widehat{\gamma}_{jn})$ in [3], $j = 2$ and 3, similarly, what results in their simple forms, and hence straightforward estimability, is the normality conditions on independent $\xi$, $\delta$, $\varepsilon$ in the SEIVM and independent $\delta$, $\varepsilon$ in the FEIVM that are required for the CLT's in [3].

**Remark 2.7.** Studentized bivariate CLT's for $(\widehat{\beta}_{jn}, \widehat{\alpha}_{jn})$, $(\widehat{\beta}_{jn}, \widehat{\gamma}_{jn})$ and $(\widehat{\alpha}_{jn}, \widehat{\gamma}_{jn})$, $1 \leq j \leq 3$, that are similar in form and features to those in (a) and (b) of Theorem 2.1 also hold true. The proofs of such CLT's are based on the auxiliary CLT's in Theorem 3.4 and are like the proofs of (a) and (b) of Theorem 2.1.

**Remark 2.8.** Condition **(B)** on the error terms in Theorem 2.1 imposes hardly any restrictions, and is only assumed there for the sake of checking (3.72) that, in particular, implies that in the (a) and (b) parts, the respective matrices $V_{z(j,n,\beta)z(j,n,\beta)}$ and $V_{z(j,n,\widehat{\beta}_{jn})z(j,n,\widehat{\beta}_{jn})}$ are positive definite on sets whose probabilities go to one, as $n \to \infty$. Moreover, along the lines of the proof of (3.72), it is not hard to see that the Studentized bivariate version of Theorem 2.1 for $(\widehat{\beta}_{jn}, \widehat{\alpha}_{jn})$ does not require assuming **(B)** at all.

**Remark 2.9.** Under (2.6) and (2.7), matrices

$$\Big(\mathrm{diag}(U(j,n), 1, L(j,n))\Big)^{-1} V_{z(j,n,\widehat{\beta}_{jn})z(j,n,\widehat{\beta}_{jn})} \Big(\mathrm{diag}(U(j,n), 1, L(j,n))\Big)^{-1}$$

are natural estimators for the respective asymptotic covariance matrices of $(\widehat{\beta}_{jn}, \widehat{\alpha}_{jn}, \widehat{\gamma}_{jn})$ obtained in [3, 7, 8] (for the proof in the case of $j = 3$ see (3.80), (3.87) and (3.89)–(3.91)). These matrices are different from the respective estimators in [3, 7, 8] that were constructed under the additional normality, or normality like, conditions specified in Remark 2.6, and also work when these conditions fail. These normality conditions also ensured positive definiteness of the asymptotic covariance matrices of $(\widehat{\beta}_{jn}, \widehat{\alpha}_{jn}, \widehat{\gamma}_{jn})$ in [3, 7, 8], $1 \leq j \leq 3$, but they seem to be more restrictive than our weak assumptions of **(B)** that guarantee such positivity on account of (3.72) and (3.89)–(3.91).

### *2.2. Interplay between* SEIVM *and* FEIVM

**Remark 2.10.** We elaborate further on assumption **(S1)** on the explanatory variables in the SEIVM (1.1). **(S1)** is inherited from the author's previous work in [14], where it was introduced the first time around for SEIVM's. In [14] that also led to [15, 17], the motivations for introducing **(S1)** into the asymptotic theory in the SEIVM (1.1) amounted to more than just aiming at a generalization of the usual assumptions in (2.6) that had been used in the literature before. From the empirical standpoint, by letting $\xi_i$ in (1.1) to have an infinite deviation ($\operatorname{Var} \xi = \infty$), we make them more dominant over the errors with finite variances. This, in turn, renders observations $y_i$ and $x_i$ to be more robust to noise (errors) and thus, more precise. Moreover, from a rigorous mathematical point of view, condition $\xi \in \mathrm{DAN}$ is optimal, or nearly optimal for the marginal CLT's for $\widehat{\beta}_{jn}$ and $\widehat{\alpha}_{jn}$ in [14, 15, 17], $1 \leq j \leq 3$ (cf., e.g., Proposition 1.1 in [14]). In additon, some distinctive features of the SEIVM's under $\operatorname{Var} \xi = \infty$ were discovered (cf., e.g., Corollary 1 and Remarks 5, 6, 7, 9 in [17]).

**Remark 2.11.** Assumptions **(F1)**–**(F3)** on the explanatory variables in the FEIVM (1.1) were first introduced and used in the FEIVM (1.1) in [15]. In [15], for the sake of achieving a strong similarity between the marginal CLT's for $\widehat{\beta}_{jn}$, $\widehat{\alpha}_{jn}$ and $\widehat{\gamma}_{jn}$ in the FEIVM (1.1) of [15] and those in the SEIVM (1.1) of [14], assumptions **(F1)**–**(F3)** on the deterministic $\xi_i$ were introduced in such a way that, due to Remarks 1.1 and 1.2, they would be natural companions for the DAN condition in **(S1)** on stochastic $\xi_i$ in [14]. An empirical rationale behind allowing $\lim_{n \to \infty} \overline{\xi^2} = \infty$ as in **(F2)** is similar to that behind possibly having $\operatorname{Var} \xi = \infty$ as in **(S1)** (cf. Remark 2.10). Also, **(F3)** is in some sense optimal for the marginal CLT's for $\widehat{\beta}_{jn}$ and $\widehat{\alpha}_{jn}$ in [15], $1 \leq j \leq 3$ (cf. Remark 2.1.8 of [15]). A further similarity of **(S1)** and **(F2)** relates to the fact that their respective partial cases $\operatorname{Var} \xi = \infty$ and $\lim_{n \to \infty} \overline{\xi^2} = \infty$ make the SEIVM (1.1) and FEIVM (1.1) behave as if they were the simple regressions $y_i = \beta x_i + \alpha + \delta_i$ (cf. Remarks 1.1.6, 2.1.10 and Sections 1.1.5, 2.1.5 in [15]).

**Remark 2.12.** Between the SEIVM under (2.6) and the FEIVM under (2.7) there is an interplay established by Gleser [8] that yields, in particular, that the CLT's for $(\widehat{\beta}_{jn}, \widehat{\alpha}_{jn}, \widehat{\gamma}_{jn})$ as in [3, 7], $j = \overline{1,3}$, that are proved in the FEIVM under (2.7) also hold true in the SEIVM with $\{\xi, \xi_i, i \geq 1\}$ satisfying (2.6). Similarly, the identity in form of the marginal CLT's for $\widehat{\beta}_{jn}$, $\widehat{\alpha}_{jn}$ and $\widehat{\gamma}_{jn}$ in the FEIVM under **(F1)**–**(F3)** in [15] to those in the SEIVM under **(S1)** of [14] establishes an asymptotic interplay between these two more general models. The CLT's for $(\widehat{\beta}_{jn}, \widehat{\alpha}_{jn}, \widehat{\gamma}_{jn})$ as in Theorem 2.1 and the respective bivariate CLT's as in Remark 2.7 are also universal in form for the latter two models, invariant as to whether the explanatory variables have a deterministic nature, as in the FEIVM, or a stochastic nature, as in the SEIVM, and thus, further contribute to the models' interplay in terms of their asymptotics.

Yu. V. Martsynyuk/CLT's in error-in-variables models 357Yu. V. Martsynyuk/CLT's in error-in-variables models  357

## 3. Auxiliary results and proofs of main results

### 3.1. CLT *for a multivariate Student statistic that is based on independent random vectors satisfying the Lindeberg condition*

In this section we state and prove Theorem 3.1, a key auxiliary CLT required for the proofs in the FEIVM of this paper. It may also be of independent interest.

For random vectors $Z_1, \cdots, Z_n$ in $\mathbb{R}^d$, we introduce a multivariate Student statistic as follows:

$$St_n(Z) = \sqrt{n}\,\overline{Z}\,V_{ZZ}^{-T/2}, \qquad (3.1)$$

where $V_{ZZ}^{1/2}$ is either the (left) Cholesky, or the symmetric positive definite, square root of the matrix $V_{ZZ} = (n-1)^{-1}\sum_{i=1}^n (Z_i - \overline{Z})^T(Z_i - \overline{Z})$. In the latter case $V_{ZZ}^{-T/2} = V_{ZZ}^{-1/2}$. Hereafter, notations $\|\cdot\|$, $\mathbb{1}_{\{\cdot\}}$ and $Z^{(j)}$ respectively stand for the Euclidean norm in $\mathbb{R}^d$, an indicator function and the $j^{\text{th}}$ component of vector $Z \in \mathbb{R}^d$, $d \geq 1$. When we will write that a (random) matrix converges (in probability) to another (random) matrix of the same size, it will mean that each entry of the converging matrix goes (in probability) to the corresponding entry of the limiting matrix.

**Theorem 3.1.** *Let $\{Z_i(n), 1 \leq i \leq n, n \geq 1\}$ be a triangular sequence of random vectors in $\mathbb{R}^d$. For each $n \geq 1$, suppose that $Z_1(n), \cdots, Z_n(n)$ are independent, $E\,Z_i(n) = 0$ and covariance matrices $\operatorname{Cov} Z_i(n)$ are finite, $1 \leq i \leq n$. Assume also that, as $n \to \infty$,*

$$\sum_{i=1}^n \operatorname{Cov} Z_i(n) \to \Sigma > 0, \qquad (3.2)$$

*with some limiting $d \times d$ matrix $\Sigma$, and that the Lindeberg condition is satisfied, namely,*

$$\text{for each } \mu > 0, \quad \sum_{i=1}^n E\left(\|Z_i(n)\|^2 \mathbb{1}_{\{\|Z_i(n)\|\geq \mu\}}\right) \to 0. \qquad (3.3)$$

*Then, for the Student statistic $St_n(Z(n))$ as in (3.1), as $n \to \infty$, $St_n(Z(n)) \xrightarrow{\mathcal{D}} N(0, I_d)$.*

*Proof.* According to the Lindeberg-Feller theorem (cf., e.g., 2.27 Proposition in [21]), (3.2) and (3.3) yield that $\sum_{i=1}^n Z_i(n) \xrightarrow{\mathcal{D}} N(0, \Sigma)$, $n \to \infty$. Since $\Sigma > 0$, Theorem 3.1 follows via Theorem 4.1 from the latter convergence in distribution and from showing that

$$(n-1)V_{Z(n)Z(n)} \xrightarrow{P} \Sigma, \quad n \to \infty. \qquad (3.4)$$

The proof of (3.4) goes first for $Z_i(n)$ that are such that the limiting matrix $\Sigma = I_d$ in (3.2).

Lindeberg's condition (3.3) implies that the same condition holds true for $\{Z_i^{(j)}(n), 1 \leq i \leq n, n \geq 1\}$ for each $1 \leq j \leq d$, that, in turn, results in

$$\frac{\sum_{i=1}^n \left(Z_i^{(j)}(n)\right)^2}{\sum_{i=1}^n \operatorname{Var} Z_i^{(j)}(n)} \xrightarrow{P} 1 \quad \text{and} \quad \frac{\sum_{i=1}^n \left(Z_i^{(j)}(n) - \overline{Z^{(j)}(n)}\right)^2}{\sum_{i=1}^n \left(Z_i^{(j)}(n)\right)^2} \xrightarrow{P} 1,$$

as $n \to \infty$ (cf., e.g., respective conclusions (3.6) and (3.7) in [16]). Consequently,

$$\frac{\sum_{i=1}^n \left(Z_i^{(j)}(n) - \overline{Z^{(j)}(n)}\right)^2}{\sum_{i=1}^n \operatorname{Var} Z_i^{(j)}(n)} \xrightarrow{P} 1, \quad n \to \infty,$$

and for having (3.4), it suffices to show the convergence in probability of the off-diagonal entries of the matrix $(n-1)V_{Z(n)Z(n)}$ to the corresponding entries of the matrix $\Sigma = I_d$, namely, for any $1 \leq j, k \leq d, j \neq k$, as $n \to \infty$,

$$\sum_{i=1}^n (Z_i^{(j)}(n) - \overline{Z^{(j)}(n)})(Z_i^{(k)}(n) - \overline{Z^{(k)}(n)}) \xrightarrow{P} 0. \tag{3.5}$$

First, for any $1 \leq j, k \leq d, j \neq k$, we prove that

$$\sum_{i=1}^n Z_i^{(j)}(n) Z_i^{(k)}(n) \xrightarrow{P} 0, \quad n \to \infty. \tag{3.6}$$

Then, for (3.6) to hold true, by Theorem 3 in [20] on p.210, it is sufficient to show that for any $\nu > 0$ and some $\tau > 0$, as $n \to \infty$,

$$\sum_{i=1}^n P\left(\left|Z_i^{(j)}(n) Z_i^{(k)}(n)\right| \geq \nu\right) \to 0, \tag{3.7}$$

$$\sum_{i=1}^n \left(E\left(\left(Z_i^{(j)}(n) Z_i^{(k)}(n)\right)^2 \mathbb{1}_{\left\{\left|Z_i^{(j)}(n) Z_i^{(k)}(n)\right| < \tau\right\}}\right)\right.$$
$$\left. - \left(E\left(Z_i^{(j)}(n) Z_i^{(k)}(n) \mathbb{1}_{\left\{\left|Z_i^{(j)}(n) Z_i^{(k)}(n)\right| < \tau\right\}}\right)\right)^2\right) \to 0 \tag{3.8}$$

and

$$\sum_{i=1}^n E\left(Z_i^{(j)}(n) Z_i^{(k)}(n) \mathbb{1}_{\left\{\left|Z_i^{(j)}(n) Z_i^{(k)}(n)\right| < \tau\right\}}\right) \to 0. \tag{3.9}$$

Since

$$\sum_{i=1}^n P\left(\left|Z_i^{(j)}(n) Z_i^{(k)}(n)\right| \geq \nu\right)$$
$$\leq \frac{1}{\nu} \sum_{i=1}^n E\left(\left|Z_i^{(j)}(n) Z_i^{(k)}(n)\right| \mathbb{1}_{\left\{\left|Z_i^{(j)}(n) Z_i^{(k)}(n)\right| \geq \nu\right\}}\right)$$
$$\leq \frac{1}{2\nu} \sum_{i=1}^n E\left(\|Z_i(n)\|^2 \mathbb{1}_{\{\|Z_i(n)\| \geq \sqrt{2\nu}\}}\right),$$

then (3.7) with any $\nu > 0$ follows from (3.3). As to (3.8), we have

$$\sum_{i=1}^{n} E\left(\left(Z_i^{(j)}(n)Z_i^{(k)}(n)\right)^2 \mathbb{1}_{\left\{\left|Z_i^{(j)}(n)Z_i^{(k)}(n)\right|<1\right\}}\right)$$

$$\leq \sum_{i=1}^{n} E\left(\left(Z_i^{(j)}(n)Z_i^{(k)}(n)\right)^2 \mathbb{1}_{\left\{\left|Z_i^{(j)}(n)\right|<1\right\}}\right)$$

$$+ \sum_{i=1}^{n} E\left(\left(Z_i^{(j)}(n)Z_i^{(k)}(n)\right)^2 \mathbb{1}_{\left\{\left|Z_i^{(k)}(n)\right|<1\right\}}\right) =: B_{1n} + B_{2n}.$$

For any $\phi > 0$ and sufficiently large $n$,

$$B_{1n} \leq \sum_{i=1}^{n} E\left(\left(Z_i^{(j)}(n)Z_i^{(k)}(n)\right)^2 \mathbb{1}_{\left\{\left|Z_i^{(j)}(n)\right|<1, \left|Z_i^{(k)}(n)\right|\geq \sqrt{\frac{\phi}{4}}\right\}}\right)$$

$$+ \sum_{i=1}^{n} E\left(\left(Z_i^{(j)}(n)Z_i^{(k)}(n)\right)^2 \mathbb{1}_{\left\{\left|Z_i^{(j)}(n)\right|<1, \left|Z_i^{(k)}(n)\right|<\sqrt{\frac{\phi}{4}}\right\}}\right)$$

$$\leq \sum_{i=1}^{n} E\left(\left(Z_i^{(k)}(n)\right)^2 \mathbb{1}_{\left\{\left|Z_i^{(k)}(n)\right|\geq \sqrt{\frac{\phi}{4}}\right\}}\right)$$

$$+ \frac{\phi}{4} \sum_{i=1}^{n} E\left(\left(Z_i^{(j)}(n)\right)^2 \mathbb{1}_{\left\{\left|Z_i^{(j)}(n)\right|<1\right\}}\right) \leq \frac{\phi}{2} + \frac{\phi}{2} = \phi,$$

where, on account of (3.2) with $\Sigma = I_d$ and that $Z_i^{(j)}(n)$ satisfy Lindeberg's condition, $\sum_{i=1}^{n} E\left((Z_i^{(j)}(n))^2 \mathbb{1}_{\{|Z_i^{(j)}(n)|<1\}}\right) \to \lim_{n\to\infty} \sum_{i=1}^{n} \operatorname{Var} Z_i^{(j)}(n) = 1$, $n \to \infty$. Consequently, $B_{1n} \to 0$, $n \to \infty$. By symmetry, $B_{2n} \to 0$, $n \to \infty$. Hence, (3.8) with $\tau = 1$ holds true. As to the convergence in (3.9), it is valid because (3.3) and the assumption that $\Sigma = I_d$ in (3.2) (hence, $\lim_{n\to\infty} \sum_{i=1}^{n} E\left(Z_i^{(j)}(n)Z_i^{(k)}(n)\right) = 0$ for any $1 \leq j, k \leq d, j \neq k$) yeild that for any $\tau > 0$,

$$\lim_{n\to\infty} \sum_{i=1}^{n} E\left(Z_i^{(j)}(n)Z_i^{(k)}(n) \mathbb{1}_{\left\{\left|Z_i^{(j)}(n)Z_i^{(k)}(n)\right|<\tau\right\}}\right)$$

$$= -\lim_{n\to\infty} \sum_{i=1}^{n} E\left(Z_i^{(j)}(n)Z_i^{(k)}(n) \mathbb{1}_{\left\{\left|Z_i^{(j)}(n)Z_i^{(k)}(n)\right|\geq\tau\right\}}\right)$$

and

$$\lim_{n\to\infty} \left|\sum_{i=1}^{n} E\left(Z_i^{(j)}(n)Z_i^{(k)}(n) \mathbb{1}_{\left\{\left|Z_i^{(j)}(n)Z_i^{(k)}(n)\right|\geq\tau\right\}}\right)\right|$$

$$\leq \lim_{n\to\infty} \sum_{i=1}^{n} E\left(\left|Z_i^{(j)}(n) Z_i^{(k)}(n)\right| \mathbb{1}_{\left\{\left|Z_i^{(j)}(n) Z_i^{(k)}(n)\right| \geq \tau\right\}}\right)$$

$$\leq \frac{1}{2} \lim_{n\to\infty} \sum_{i=1}^{n} E\left(\|Z_i(n)\|^2 \mathbb{1}_{\{\|Z_i(n)\| \geq \sqrt{2\tau}\}}\right) = 0.$$

This completes the proof of (3.6).

Now, since (3.6) holds true and $\sum_{i=1}^{n}\left(Z_i^{(j)}(n) - \overline{Z^{(j)}(n)}\right)\left(Z_i^{(k)}(n) - \overline{Z^{(k)}(n)}\right)$ $= \sum_{i=1}^{n} Z_i^{(j)}(n) Z_i^{(k)}(n) - n\,\overline{Z^{(j)}(n)}\,\overline{Z^{(k)}(n)}$, to show (3.5), it suffices to prove that for any $1 \leq j, k \leq d$, $j \neq k$,

$$n\,\overline{Z^{(j)}(n)}\,\overline{Z^{(k)}(n)} \xrightarrow{P} 0, \quad n \to \infty.$$

The latter follows from (3.2) and Markov's inequality for any $1 \leq j \leq d$ and $\phi > 0$, namely,

$$P\left(n\left(\overline{Z^{(j)}(n)}\right)^2 \geq \phi \sum_{i=1}^{n} \operatorname{Var} Z_i^{(j)}(n)\right) \leq (\phi n)^{-1} \to 0, \quad n \to \infty.$$

This also completes the proof of (3.4) for $Z_i(n)$ with $\Sigma = I_d$ in (3.2).

When $\Sigma \neq I_d$ in (3.2), we take $\{\widetilde{Z}_i(n) = Z_i(n)\Sigma^{-T/2}, 1 \leq i \leq n, n \geq 1\}$. Clearly, $\widetilde{Z}_i(n)$ satisfy all the conditions of Theorem 3.1, with respective $\Sigma = I_d$ in (3.2). Hence, (3.4) with $\widetilde{Z}_i(n)$ replacing $Z_i(n)$ holds true, namely, $(n-1)V_{\widetilde{Z}(n)\widetilde{Z}(n)} \xrightarrow{P} I_d$, $n \to \infty$. The latter convergence implies (3.4) for $Z_i(n)$. □

**Remark 3.1.** Though the CLT of Theorem 3.1 appears to be quite natural, especially in view of the well-known multivariate Lindeberg-Feller CLT (cf., e.g., 2.27 Proposition in [21]) and Theorem 4.1, one essential link for its conclusion, namely, an appropriate version of (4.1) as in (3.4) used to be missing. Our (3.4) is essentially a multivariate extension of a part of Raikov's theorem (cf. Theorem 4 on p.143 in [9]) that amounts to saying that for r.v.'s $\{X_i(n), 1 \leq i \leq n, n \geq 1\}$ that are independent in each row and satisfy the Lindeberg condition, we have $\sum_{i=1}^{n}(X_i(n) - E\,X_i(n))^2 \Big/ \sum_{i=1}^{n} \operatorname{Var} X_i(n) \xrightarrow{P} 1$, $n \to \infty$. On p.145 in [22], the authors of Theorem 4.1 as in Section 4 of the present paper, with applications of their theorem in mind, pose the question of having (4.1) with some matrix $V_n$ as in Theorem 4.1 that would correspond to a case of independent nonidentically distributed random vectors, at least when the latter have finite covariance matrices. Hence, our (3.4) can also be viewed as a partial answer to this question.

### 3.2. *Special characterization of the generalized domain of attraction of the multivariate normal law*

The purpose of this subsection is to establish a special, convenient characterization of the generalized domain of attraction of the multivariate normal law (GDAN) as in Theorem 3.2. It will enable us to apply Theorem 4.3 for obtaining an auxiliary CLT of (a) of Theorem 3.3 for proving (a) of Theorem 2.1 in the SEIVM (1.1). By establishing Theorem 3.2, we also give an example of special vectors $Z = (Z^{(1)}, \cdots, Z^{(d)})$ for which the fact that $Z^{(j)} \in \text{DAN}$ for all $1 \leq j \leq d$ characterizes that $Z \in \text{GDAN}$. In this example we get away from the condition that components $Z^{(j)}$ of vector $Z$ are identically distributed as in the related example with spherically symmetric $Z$ at the end of Remark 4.5.

We recall that random vector $Z$ is a full vector, or has a full distribution, if $\langle Z, u \rangle$ is a nondegenerate r.v. for all deterministic unit norm vectors $u$. In the sequel, for $Z \in \mathbb{R}^d$, when $d \geq 2$, $Z^{(k,k+l)} = (Z^{(k)}, Z^{(k+1)}, \cdots, Z^{(k+l)})$ denotes a subvector of $Z$ that has all the components of $Z$ starting with $Z^{(k)}$ and ending with $Z^{(k+l)}$, $1 \leq k \leq d-1$, $1 \leq l \leq k+l \leq d$. Notation $\text{diag}(\cdot, \cdots, \cdot)$ stands for a block-diagonal matrix, where in the brackets square matrix blocks that are on its diagonal are listed.

**Theorem 3.2.** *Let $\{Z, Z_i, i \geq 1\}$ be i.i.d. random vectors in $\mathbb{R}^d$. Suppose that*

$$\text{for } j \neq k, \quad E|Z^{(j)}Z^{(k)}| < \infty, \quad \text{if} \quad E(Z^{(j)})^2 = \infty \quad \text{and/or} \quad E(Z^{(k)})^2 = \infty. \tag{3.10}$$

*For vector $\widetilde{Z}$ formed by all the components, if any, of $Z$ whose second moments exist, assume that*

$$\widetilde{Z} \text{ is full.} \tag{3.11}$$

*Then the following two statements are equivalent:*

(a) $Z \in \text{GDAN}$;
(b) $Z^{(j)} \in \text{DAN}$ *for all* $1 \leq j \leq d$.

*Proof.* That (a) implies (b) is explained in Remark 4.5.

Conversely, assume that (b) holds true. By equivalence of the (a) and (b) parts of Theorem 4.2, the proof of (a) of this lemma reduces to verifying convergence in (b) of Theorem 4.2 for suitably chosen matrices $B_n$, provided that vector $Z$ is full. Clearly, we are concerned with the case of $d \geq 2$ only.

If $E(Z^{(j)})^2 < \infty$ for all $1 \leq j \leq d$, then $Z = \widetilde{Z}$ is full by (3.11). Hence, $\text{Cov } Z > 0$ and, due to a weak law of large numbers applied to each entry of $\sum_{i=1}^n (Z_i - \overline{Z})^T(Z_i - \overline{Z})$, the (b) part of Theorem 4.2 is satisfied with matrix

$$B_n = n^{-1/2} \left( \text{Cov } Z \right)^{-1/2}. \tag{3.12}$$

Suppose now that, without loss of generality, $E(Z^{(j)})^2 = \infty$ for all $1 \leq j \leq m$, $1 \leq m < d$, while $E(Z^{(j)})^2 < \infty$ for all $m+1 \leq j \leq d$. First, note that such vector $Z$ is full and thus, Theorem 4.2 is applicable. Indeed, for any unit norm scalar vector $u$, $\text{Var}\langle Z, u \rangle = \sum_{j=1}^d (u^{(j)})^2 \text{Var } Z^{(j)} + \sum_{j \neq k} u^{(j)} u^{(k)} \text{cov}(Z^{(j)}, Z^{(k)})$, and, if $|u^{(1)}| + \cdots + |u^{(m)}| > 0$, on account of (3.10), $\text{Var}\langle Z, u \rangle = \infty$, while when $u^{(1)} = \cdots = u^{(m)} = 0$, $\text{Var}\langle Z, u \rangle > 0$ by (3.11). Next, we introduce the block-diagonal matrix

$$B_n = n^{-1/2} \text{diag}\Big((\ell^{(1)}(n))^{-1}, \cdots, (\ell^{(m)}(n))^{-1}, \left(\text{Cov } Z^{(m+1,\,d)}\right)^{-1/2}\Big), \tag{3.13}$$

where $\ell^{(1)}(n) \nearrow \infty, \cdots, \ell^{(m)}(n) \nearrow \infty$ are slowly varying functions at infinity that correspond to $Z^{(j)} \in$ DAN (cf. Remark 1.1), namely, we have for $1 \leq j \leq m$,

$$\frac{\sum_{i=1}^n (Z_i^{(j)} - EZ^{(j)})}{\sqrt{n}\ell^{(j)}(n)} \xrightarrow{\mathcal{D}} N(0,1), \quad n \to \infty. \tag{3.14}$$

Note that $\left(\operatorname{Cov} Z^{(m+1,d)}\right)^{-1/2}$ in (3.13) is well-defined on account of (3.11). For $B_n$ as in (3.13), we are to verify convergence in the (b) part of Theorem 4.2, namely, that

$$B_n \sum_{i=1}^n (Z_i - \overline{Z})^T (Z_i - \overline{Z}) B_n^T = B_n \sum_{i=1}^n s_{i,ZZ} B_n^T$$

$$:= E_n = (e_n^{jk})_{j,k=\overline{1,d}} = \begin{pmatrix} E_n^1 & E_n^2 \\ (E_n^2)^T & E_n^3 \end{pmatrix} \xrightarrow{P} I_d, \quad n \to \infty, \tag{3.15}$$

where matrix $E_n$ consists of the following matrix blocks:

$$E_n^1 = (e_n^{jk})_{1 \leq j,k \leq m}, \quad \text{with} \quad e_n^{jk} = \frac{\sum_{i=1}^n s_{i,Z^{(j)}Z^{(k)}}}{n\ell^{(j)}(n)\ell^{(k)}(n)}, \tag{3.16}$$

$$E_n^2 = (e_n^{jk})_{1 \leq j \leq m,\, m+1 \leq k \leq d}, \quad e_n^{jk} = \operatorname{const} \frac{\sum_{i=1}^n s_{i,Z^{(j)}Z^{(m+1)}}}{n\ell^{(j)}(n)}$$

$$+ \cdots + \operatorname{const} \frac{\sum_{i=1}^n s_{i,Z^{(j)}Z^{(d)}}}{n\ell^{(j)}(n)}, \tag{3.17}$$

with some absolute constants in (3.17) depending on $j$ and $k$, but not on $n$, and

$$E_n^3 = \left(\operatorname{Cov} Z^{(m+1,d)}\right)^{-\frac{1}{2}} \frac{\sum_{i=1}^n s_{i,Z^{(m+1,d)}Z^{(m+1,d)}}}{n} \left(\operatorname{Cov} Z^{(m+1,d)}\right)^{-\frac{T}{2}} \tag{3.18}$$

As $n \to \infty$, on account of (3.14) and (4.7) of Remark 4.4,

$$e_n^{jj} \xrightarrow{P} 1 \quad \text{for all} \quad 1 \leq j \leq m, \tag{3.19}$$

and, due to (3.10) and the fact that $\ell^{(j)}(n) \nearrow \infty$, $1 \leq j \leq m$,

$$e_n^{jk} \xrightarrow{P} 0 \quad \text{for} \quad j \neq k,\ 1 \leq j \leq m \text{ and } 1 \leq k \leq d, \tag{3.20}$$

while, clearly,

$$E_n^3 \xrightarrow{P} I_{d-m}. \tag{3.21}$$

This concludes the proof of the convergence in (3.15) with $B_n$ as in (3.13).

Finally, the third type of vectors satisfying (3.10), (3.11) and having all their components in DAN consists of vectors $Z$ with $E(Z^{(j)})^2 = \infty$ for all $1 \leq j \leq d$. Such vectors $Z$ are full, since for any deterministic unit norm vector $u$, on account of (3.10), $\operatorname{Var}\langle Z, u \rangle = \infty$. Put

$$B_n = n^{-1/2} \operatorname{diag}\left((\ell^{(1)}(n))^{-1}, \cdots, (\ell^{(d)}(n))^{-1}\right), \tag{3.22}$$

where $\ell^{(j)}(n)$ are as in (3.14), and $\ell^{(j)}(n) \nearrow \infty$, $n \to \infty$, for all $1 \leq j \leq d$. For all elements $e^{jk}$ of matrix $E_n$ in (3.15) defined now via $B_n$ of (3.22), by (4.7) and (3.10),

$$e_n^{jk} = \frac{\sum_{i=1}^n s_{i,Z^{(j)}Z^{(k)}}}{n\ell^{(j)}(n)\ell^{(k)}(n)} \xrightarrow{P} \begin{cases} 1, \, j = k, \\ 0, \, j \neq k, \end{cases} \quad 1 \leq j, k \leq d, \qquad (3.23)$$

as $n \to \infty$. Hence, convergence in probability in (3.15) with $B_n$ of (3.22) is valid. $\square$

**Remark 3.2.** Relationship in (3.10) may follow from independence of $Z^{(j)}$ and $Z^{(k)}$ and the existence of their respective first moments. According to the Cauchy-Schwarz inequality, (3.10) is also satisfied when, e.g., $Z^{(1)} \in$ DAN with $E(Z^{(1)})^2 = \infty$ (from Remark 1.1, $E(Z^{(1)})^{2-a} < \infty$ for any $a \in (0,2]$), while $E(Z^{(j)})^{2+\Delta} < \infty$ for some $\Delta > 0$, for all $2 \leq j \leq d$.

### 3.3. Auxiliary CLT's required for the proof of Theorem 2.1 in Section 3.4

This subsection presents Theorems 3.3 and 3.4 containing auxiliary CLT's required for the proof of Theorem 2.1 in Section 3.4. Theorems 3.3 and 3.4 are based on, and extend and unite, the auxiliary results in [14, 15, 16, 17]. This subsection is written with the genuinely multivariate case of $d > 1$ in mind. This is understood from, but not spelled out in, some of the conditions and notations that are to appear (cf., e.g., (3.25)). However, simply omitting the arguments that are suitable and used only for the case of $d > 1$, makes the results of this subsection also valid for the case of $d = 1$.

Let

$$\zeta_i = \Big((\xi_i - m)\delta_i, (\xi_i - m)\varepsilon_i, \delta_i, \varepsilon_i, \delta_i\varepsilon_i - \mu, \delta_i^2 - \lambda\theta, \varepsilon_i^2 - \theta\Big), \quad 1 \leq i \leq n, \quad (3.24)$$

where

$$m = \begin{cases} E\xi, & \text{in the SEIVM (1.1)}, \\ \lim_{n\to\infty} \overline{\xi}, & \text{in the FEIVM (1.1)}. \end{cases}$$

Using (3.24) and vectors of constants $b_1, \cdots, b_d$ in $\mathbb{R}^7$, whose components are such that

$$\text{if } \begin{cases} \operatorname{Var}\xi = \infty, & \text{in the SEIVM}, \\ \lim_{n\to\infty} \overline{\xi^2} = \infty, & \text{in the FEIVM}, \end{cases} \text{ then } b_j^{(1)} = b_j^{(2)} = 0 \text{ for all } 2 \leq j \leq d, \tag{3.25}$$

we define vectors

$$J_i = \Big(\langle \zeta_i, b_1 \rangle, \cdots, \langle \zeta_i, b_d \rangle\Big), \quad 1 \leq i \leq n. \qquad (3.26)$$

In the SEIVM (1.1), we will also consider random vectors $\zeta_0$ and $J_0$ that respectively generate $\zeta_i$ and $J_i$, namely,

$$\zeta_0 = \Big((\xi - m)\delta, (\xi - m)\varepsilon, \delta, \varepsilon, \delta\varepsilon - \mu, \delta^2 - \lambda\theta, \varepsilon^2 - \theta\Big) \qquad (3.27)$$

and
$$J_0 = \Big(\langle \zeta_0, b_1\rangle, \cdots, \langle \zeta_0, b_d\rangle\Big). \tag{3.28}$$

Next, for a special case of the multivariate Student statistic in (3.1), we prove Theorem 3.3, the first of the two auxiliary theorems of this subsection.

**Theorem 3.3.** *Consider* $St_n(J) = \sqrt{n}\,\overline{J}V_{JJ}^{-T/2}$ *based on* $J_i$ *of (3.26) that are defined via nonrandom vectors* $b_1, \cdots, b_d$ *in* $\mathbb{R}^7$ *satisfying (3.25).*

(a) *In the SEIVM (1.1), assume* **(A)**, **(S1)**, **(S2)** *and also that*
$$\begin{cases} J_0 \text{ is full}, & \text{if } \mathrm{Var}\,\xi < \infty \text{ and/or } b_1^{(1)} = b_1^{(2)} = 0, \\ J_0^{(2,d)} = \Big(\langle \zeta_0, b_2\rangle, \cdots, \langle \zeta_0, b_d\rangle\Big) \text{ is full}, & \text{if } \mathrm{Var}\,\xi = \infty \text{ and } |b_1^{(1)}| + |b_1^{(2)}| > 0. \end{cases} \tag{3.29}$$

*Then, as* $n \to \infty$, $St_n(J) \xrightarrow{\mathcal{D}} N(0, I_d)$.

(b) *In the FEIVM (1.1), let* **(A)** *and* **(F1)**–**(F3)** *be satisfied, and*
$$\begin{cases} \lim_{n\to\infty} \dfrac{1}{n}\sum_{i=1}^{n} \mathrm{Cov}\,J_i > 0, & \text{if } \lim_{n\to\infty} \overline{\xi^2} < \infty \text{ and/or } b_1^{(1)} = b_1^{(2)} = 0, \\ \lim_{n\to\infty} \dfrac{1}{n}\sum_{i=1}^{n} \mathrm{Cov}\,J_i^{(2,d)} > 0, & \text{if } \lim_{n\to\infty} \overline{\xi^2} = \infty \text{ and } |b_1^{(1)}| + |b_1^{(2)}| > 0. \end{cases} \tag{3.30}$$

*Then, as* $n \to \infty$, $St_n(J) \xrightarrow{\mathcal{D}} N(0, I_d)$.

*Proof of the* (a) *part of Theorem* 3.3. The proof follows from Theorem 4.3, provided that we show that $J_0 \in \mathrm{GDAN}$ by using Theorem 3.2.

First, we argue that the components of $J_0$ obey (3.10) and (3.11). Indeed, if $\mathrm{Var}\,\xi < \infty$ and/or $b_1^{(1)} = b_1^{(2)} = 0$, then $E\langle \zeta_0, b_j\rangle^2 < \infty$ for all $1 \leq j \leq d$, and due to (3.29), conditions (3.10) and (3.11) are clear satisfied. If $\mathrm{Var}\,\xi = \infty$ and $|b_1^{(1)}|+|b_1^{(2)}| > 0$, then (3.25) implies that $E\langle \zeta_0, b_1\rangle^2 = \infty$, while $E\langle \zeta_0, b_j\rangle^2 < \infty$ for all $2 \leq j \leq d$. In this case (3.11) is guaranteed by (3.29), while (3.10) follows from $E\big|\langle \zeta_0, b_1\rangle\langle \zeta_0, b_j\rangle\big| < \infty$ for all $2 \leq j \leq d$, that can be checked straightforwardly.

Further, from (26) of the proof of Lemma 6 in [17], by noticing that condition (24) with $\langle \zeta_0, b_j\rangle$ in place of $\langle \zeta_0, b\rangle$ in that lemma is now a part of (3.29), we have that $\langle \zeta_0, b_j\rangle \in \mathrm{DAN}$, for all $1 \leq j \leq d$. Finally, on account of (b) implying (a) in Theorem 3.2, $J_0 \in \mathrm{GDAN}$. $\square$

*Proof of the* (b) *part of Theorem* 3.3. Consider a triangular sequence of random vectors $\{Z_i(n), 1 \leq i \leq n, n \geq 1\}$, where

$Z_i(n) =$
$$\begin{cases} \dfrac{J_i}{\sqrt{n}}, & \text{if } \lim_{n\to\infty} \overline{\xi^2} < \infty \text{ and/or } b_1^{(1)} = b_1^{(2)} = 0, \\ \left(\dfrac{J_i^{(1)}}{\sqrt{\sum_{i=1}^{n}(\xi_i - m)^2}}, \dfrac{J_i^{(2,d)}}{\sqrt{n}}\right), & \text{if } \lim_{n\to\infty} \overline{\xi^2} = \infty \text{ and } |b_1^{(1)}| + |b_1^{(2)}| > 0, \end{cases}$$

(3.31)

and $Z_1(n), \cdots, Z_n(n)$ are independent with $E Z_i(n) = 0$ and finite Cov $Z_i(n)$.

For $Z_i(n)$ of (3.31), we are to check conditions (3.2) and (3.3) of Theorem 3.1. If $\lim_{n\to\infty} \overline{\xi^2} < \infty$ and/or $b_1^{(1)} = b_1^{(2)} = 0$, then by (3.30), $\lim_{n\to\infty} \sum_{i=1}^n \text{Cov } Z_i(n) = \lim_{n\to\infty} \sum_{i=1}^n \text{Cov } J_i/n > 0$. Otherwise, it is not hard to see that $\lim_{n\to\infty} \sum_{i=1}^n \text{Cov } Z_i(n) = \text{diag}\left(\langle b_1^{(1,2)}\Gamma, b_1^{(1,2)}\rangle, \lim_{n\to\infty} n^{-1} \sum_{i=1}^n \text{Cov } J_i^{(2,d)}\right)$, since by (3.25), $\lim_{n\to\infty} \sum_{i=1}^n \text{cov}(Z_i^{(1)}(n), Z_i^{(j)}(n)) = \lim_{n\to\infty}(\text{const} \sum_{i=1}^n (\xi_i - m) + \text{const} \cdot n)(n \sum_{i=1}^n (\xi_i - m)^2)^{-1/2} = 0$ in this case. Hence, by **(A)** and (3.30), (3.2) is also satisfied when $\lim_{n\to\infty} \overline{\xi^2} = \infty$ and $|b_1^{(1)}| + |b_1^{(2)}| > 0$. As to the Lindeberg condition (3.3) for (3.31), it suffices to show that for all $1 \le j, k \le d$, as $n \to \infty$,

$$\text{for any } \mu > 0, \quad \sum_{i=1}^n E\left((Z_i^{(j)}(n))^2 \mathbb{1}_{\{(Z_i^{(k)}(n))^2 \ge \mu\}}\right) \to 0. \tag{3.32}$$

Let

$$\gamma_{in} = \left(\frac{\zeta_i^{(1,2)}}{\sqrt{\sum_{i=1}^n (\xi_i - m)^2}}, \frac{\zeta_i^{(3,7)}}{\sqrt{n}}\right) \tag{3.33}$$

and

$$\widetilde{b}_{jn} = \left(\sqrt{\sum_{i=1}^n (\xi_i - m)^2}\, b_j^{(1,2)}, \sqrt{n}\, b_j^{(3,7)}\right), \quad 1 \le j \le d. \tag{3.34}$$

Then,

$$\langle \zeta_i, b_j \rangle = \langle \gamma_{in}, \widetilde{b}_{jn} \rangle. \tag{3.35}$$

If $\lim_{n\to\infty} \overline{\xi^2} < \infty$ and/or $b_1^{(1)} = b_1^{(2)} = 0$, then by (3.35),

$$\sum_{i=1}^n E\left((Z_i^{(j)}(n))^2 \mathbb{1}_{\{(Z_i^{(k)}(n))^2 \ge \mu\}}\right) = n^{-1} \sum_{i=1}^n E\left(\langle \zeta_i, b_j \rangle^2 \mathbb{1}_{\{\langle \zeta_i, b_k \rangle^2 \ge \mu n\}}\right)$$

$$\le n^{-1} \sum_{i=1}^n E\left(\|\gamma_{in}\|^2 \|\widetilde{b}_{jn}\|^2 \mathbb{1}_{\{\|\gamma_{in}\|^2 \|\widetilde{b}_{kn}\|^2 \ge \mu n\}}\right)$$

$$\le \text{const} \sum_{i=1}^n E\left(\|\gamma_{in}\|^2 \mathbb{1}_{\{\text{const } \|\gamma_{in}\|^2 \ge \mu\}}\right),$$

where for any $1 \le j \le d$ and large enough $n$,

$$\frac{\|\widetilde{b}_{jn}\|^2}{n} = \frac{\text{const} \sum_{i=1}^n (\xi_i - m)^2 + \text{const} \cdot n}{n} \le \text{const}.$$

Consequently, (3.32) reduces to the Lindeberg condition for $\{\|\gamma_{in}\|, 1 \le i \le n, n \ge 1\}$ that amounts to having the verified (3.36) of the proof of Lemma 3.4 in [16] (we note that condition **(F2)** is replaced there by a weaker one, namely, $\liminf_{n\to\infty}(\overline{\xi^2} - (\overline{\xi})^2) > 0$). In the case of $\lim_{n\to\infty}\overline{\xi^2} = \infty$ and $|b_1^{(1)}| + |b_1^{(2)}| > 0$, the arguments for (3.32) are similar, and also lead to the just mentioned Lindeberg condition by additionally using (3.25).

Thus, conditions (3.2) and (3.3) for (3.31) are satisfied, and by Theorem 3.1,

$$St_n(Z(n)) \xrightarrow{\mathcal{D}} N(0, I_d), \quad n \to \infty. \tag{3.36}$$

Now, if $\lim_{n\to\infty}\overline{\xi^2} < \infty$ and/or $b_1^{(1)} = b_1^{(2)} = 0$, then $St_n(Z(n)) = St_n(J)$ and (3.36) implies $St_n(J) \xrightarrow{\mathcal{D}} N(0, I_d)$, $n \to \infty$. Otherwise, concluding the latter convergence from (3.36) requires some work that takes the rest of the proof.

Assume that $\lim_{n\to\infty}\overline{\xi^2} = \infty$ and $|b_1^{(1)}| + |b_1^{(2)}| > 0$. Define matrix

$$D_n = \mathrm{diag}\left(\frac{1}{\sqrt{\sum_{i=1}^n (\xi_i - m)^2}}, \frac{I_{d-1}}{\sqrt{n}}\right). \tag{3.37}$$

Then $Z_i(n) = J_i D_n$ and

$$St_n(Z(n)) = \sqrt{n}\,\overline{J} D_n V_{Z(n)Z(n)}^{-T/2} = St_n(J) V_{JJ}^{T/2} D_n V_{Z(n)Z(n)}^{-T/2}, \tag{3.38}$$

where matrix $V_{Z(n)Z(n)} = (n-1)^{-1}\sum_{i=1}^n (Z_i(n) - \overline{Z(n)})^T (Z_i(n) - \overline{Z(n)})$. It was shown earlier in the proof that $Z_i(n)$ of (3.31) obey the conditions of Theorem 3.1. Combining (3.4) and (4.5) with $A_n = \Sigma > 0$ therein, where matrix $\Sigma$ is as in (3.2), we conclude that

$$\left((n-1)V_{Z(n)Z(n)}\right)^{-1/2} \xrightarrow{P} \Sigma^{-1/2},$$

and hence, also, that

$$\left((n-1)V_{Z(n)Z(n)}\right)^{T/2} = \left((n-1)V_{Z(n)Z(n)}\right)^{-1/2}\left((n-1)V_{Z(n)Z(n)}\right) \xrightarrow{P} \Sigma^{T/2}, \tag{3.39}$$

as $n \to \infty$. Next, we observe that in view of (3.4),

$$\Sigma^{-1/2} D_n (n-1) V_{JJ} D_n \Sigma^{-T/2} = \Sigma^{-1/2}(n-1)V_{Z(n)Z(n)}\Sigma^{-T/2} \xrightarrow{P} I_d, \quad n \to \infty. \tag{3.40}$$

Now, (3.40) and (4.5) with $A_n^{1/2} = D_n^{-1}\Sigma^{1/2}$ ($A_n = D_n^{-1}\Sigma D_n^{-1} > 0$) therein result in, as $n \to \infty$,

$$\Sigma^{T/2} D_n^{-1}\left((n-1)V_{JJ}\right)^{-T/2} \xrightarrow{P} I_d,$$

or, equivalently, in

$$D_n^{-1}\left((n-1)V_{JJ}\right)^{-T/2} \xrightarrow{P} \Sigma^{-T/2}. \tag{3.41}$$

From (3.39) and (3.41), as $n \to \infty$,

$$\left(V_{JJ}^{T/2} D_n V_{Z(n)Z(n)}^{-T/2}\right)^{-1} = \left((n-1)V_{Z(n)Z(n)}\right)^{T/2} D_n^{-1}\left((n-1)V_{JJ}\right)^{-T/2} \xrightarrow{P} I_d,$$

while the latter convergence, (3.38) and (3.36) yield $St_n(J) \xrightarrow{\mathcal{D}} N(0, I_d)$, as $n \to \infty$. $\square$

**Remark 3.3.** In general, condition (3.29) is checked on a case-by-case basis, depending on the vectors $b_1, \cdots, b_d$ in hand. The first line of (3.29) amounts to saying that r.v. $\langle \zeta_0, u_1 b_1 + \cdots + u_d b_d \rangle$ is nondegenerate for any vector $u = (u_1, \cdots, u_d)$, $\|u\| = 1$. In particular, on account of the proof of Lemma 6 in [17], the latter statement holds true when $\operatorname{Var} \xi < \infty$ and the first two components of vector $u_1 b_1 + \cdots + u_d b_d$ are not simultaneously zero for any $u = (u_1, \cdots, u_d)$, $\|u\| = 1$, namely, when the nonrandom vectors $b_1^{(1,2)}, \cdots, b_d^{(1,2)}$ are linearly independent. The assumption in the second line of (3.29) is interpreted similarly. As to the conditions in (3.30), due to the similarity in form of the matrices $\operatorname{Cov} J_0$ and $\operatorname{Cov} J_0^{(2,d)}$ in the SEIVM respectively to those of $\lim_{n \to \infty} n^{-1} \sum_{i=1}^n \operatorname{Cov} J_i$ and $\lim_{n \to \infty} n^{-1} \sum_{i=1}^n \operatorname{Cov} J_i^{(2,d)}$ in the FEIVM, (3.30) holds true whenever (3.29) does. Indeed, the respective entries of the aforementioned corresponding matrices are the same functions of the moments and cross-moments of order $\leq 4$ of the error terms, and of $m = E\xi$ and, in case of $\operatorname{Cov} J_0$, also of $M = E\xi^2$, when dealing with the SEIVM (1.1), and, respectively, of $m = \lim_{n \to \infty} \overline{\xi}$ and, in case of $\lim_{n \to \infty} n^{-1} \sum_{i=1}^n \operatorname{Cov} J_i$, also of $M = \lim_{n \to \infty} \overline{\xi^2}$ when the FEIVM (1.1) obtains. Moreover, the assumptions in the FEIVM (1.1) and SEIVM (1.1) on the error moments, $m$ and $M$ are in complete synchrony (cf. **(A)**, **(S1)**, **(F1)**, **(F2)**).

On using Theorem 3.3, we are to study another Studentized partial sum in Theorem 3.4, the second auxiliary theorem of this subsection. This Studentized partial sum is a prototype for the main terms in the expansions for the Studentized estimators of $(\beta, \alpha, \gamma)$ as in Theorem 2.1.

Let

$$\eta_i(n) = (y_i - \alpha, x_i, s_{i,yy} - \lambda\theta, s_{i,xy} - \mu, s_{i,xx} - \theta), \quad 1 \leq i \leq n, \qquad (3.42)$$

$c_1, c_2, \cdots, c_d \in \mathbb{R}^5$ be nonzero vectors of constants, and

$$K_i(n) = \Big( \langle \eta_i(n), c_1 \rangle, \cdots, \langle \eta_i(n), c_d \rangle \Big), \quad 1 \leq i \leq n. \qquad (3.43)$$

**Theorem 3.4.** *Consider the random vectors in (3.43) defined via (3.42) and nonzero deterministic vectors $c_1, \cdots, c_d$ in $\mathbb{R}^5$ satisfying*

$$c_j^{(1)}\beta + c_j^{(2)} = 0 \quad \text{and} \quad c_j^{(3)}\beta^2 + c_j^{(4)}\beta + c_j^{(5)} = 0 \quad \text{for all } 1 \leq j \leq d. \quad (3.44)$$

*Assume also (3.25) for vectors $b_1, \cdots, b_d$ in $\mathbb{R}^7$ defined by*

$$b_j = \Big( 2\beta c_j^{(3)} + c_j^{(4)}, \, \beta c_j^{(4)} + 2c_j^{(5)}, \, c_j^{(1)}, \, c_j^{(2)}, \, c_j^{(4)}, \, c_j^{(3)}, \, c_j^{(5)} \Big), \quad 1 \leq j \leq d. \quad (3.45)$$

(a) *Let **(A)**, **(S1)** and **(S2)** be valid in the SEIVM (1.1). Suppose also that (3.29) holds true. Then for $St_n(K(n)) = \sqrt{n}\, \overline{K(n)} V_{K(n)K(n)}^{-T/2} = \sqrt{n}\, \overline{K(n)}$ $\Big( (n-1)^{-1} \sum_{i=1}^n (K_i(n) - \overline{K(n)})^T (K_i(n) - \overline{K(n)}) \Big)^{-T/2}$, as $n \to \infty$, $St_n(K(n)) \xrightarrow{\mathcal{D}} N(0, I_d)$.*

(b) *Suppose that conditions* **(A)**, **(F1)**–**(F3)** *and (3.30) in the* FEIVM *(1.1) are satisfied. Then, as* $n \to \infty$, $St_n(K(n)) \xrightarrow{\mathcal{D}} N(0, I_d)$.

*Proof of the* (a) *part of Theorem* 3.4. We have,

$$St_n(K(n)) = (St_n(J) + R(n)V_{JJ}^{-T/2})V_{JJ}^{T/2}V_{K(n)K(n)}^{-T/2}, \quad (3.46)$$

where vector $R(n) = \sqrt{n}\,\overline{K(n)} - \sqrt{n}\,\overline{J}$, and, according to (41) of [17] and (3.25) for $b_j$, components $R^{(j)}(n)$ of vector $R(n)$ are such that

$$\begin{cases} R^{(1)}(n) = o_P(1), & \text{if } \operatorname{Var}\xi < \infty \text{ and/or } b_1^{(1)} = b_1^{(2)} = 0, \\ \dfrac{R^{(1)}(n)}{\ell_\xi(n)} = o_P(1), & \text{if } \operatorname{Var}\xi = \infty \text{ and } |b_1^{(1)}| + |b_1^{(2)}| > 0, \\ R^{(j)}(n) = o_P(1), & \text{for all } 2 \leq j \leq d, \end{cases} \quad (3.47)$$

with slowly varying function $\ell_\xi(n)$ as in Remark 1.1, such that $\ell_\xi(n) \nearrow \infty$ when $\operatorname{Var}\xi = \infty$, $n \to \infty$. If in (3.46), as $n \to \infty$,

$$\left\|R(n)V_{JJ}^{-T/2}\right\| = o_P(1) \quad (3.48)$$

and

$$V_{JJ}^{T/2}V_{K(n)K(n)}^{-T/2} \xrightarrow{P} I_d, \quad (3.49)$$

then the (a) part of Theorem 3.4 follows from that of Theorem 3.3 for $St_n(J)$.

First, we are to show (3.48). Introduce matrix

$$B_n = \begin{cases} n^{-\frac{1}{2}}\operatorname{diag}\left(\left(\ell_\xi(n)\sqrt{\operatorname{Var}(b_1^{(1)}\delta + b_1^{(2)}\varepsilon)}\right)^{-1}, (\operatorname{Cov}J_0^{(2,\,d)})^{-\frac{1}{2}}\right), \\ \hspace{4cm} \text{if } \operatorname{Var}\xi = \infty, |b_1^{(1)}| + |b_1^{(2)}| > 0, \\ n^{-\frac{1}{2}}(\operatorname{Cov}J_0)^{-\frac{1}{2}} \hspace{2cm}, \text{otherwise}. \end{cases} \quad (3.50)$$

Note that $B_n$ is well-defined on account of **(A)** and (3.29). Interpreting (3.48) as a degenerate weak convergence in $(\mathbb{R}^d, \|\cdot\|)$, by Theorem 4.1 from the Appendix, for (3.48) it suffices to show that, as $n \to \infty$,

$$B_n(n-1)V_{JJ}B_n^T \xrightarrow{P} I_d \quad (3.51)$$

and

$$\|\sqrt{n-1}R(n)B_n^T\| = o_P(1). \quad (3.52)$$

Convergence in (3.51) follows from the fact that $J_0$ is as in (3.10) and (3.11) (this was shown in the proof of the (a) part of Theorem 3.3), and is argued the same way as convergence in (3.15) with therein matrices $B_n$ as in (3.12) or (3.13). In this regard, we also note that the correspondence of $B_n$ in (3.50) to $B_n$ in (3.12) or (3.13) is seen by noticing that if $\operatorname{Var}\xi = \infty$ and $|b_1^{(1)}| + |b_1^{(2)}| > 0$, then, on account of Remark 10 in [17], $\sum_{i=1}^n J_i^{(1)}\left(\sqrt{n\operatorname{Var}(b_1^{(1)}\delta + b_1^{(2)}\varepsilon)}\ell_\xi(n)\right)^{-1} \xrightarrow{\mathcal{D}}$ $N(0, 1)$, $n \to \infty$. As to (3.52), when $\operatorname{Var} \xi < \infty$ and/or $b_1^{(1)} = b_1^{(2)} = 0$, it is a direct consequence of (3.47). Otherwise, namely under $\operatorname{Var} \xi = \infty$ and $|b_1^{(1)}| + |b_1^{(2)}| > 0$, it is due to (3.47) and observing that

$$\left\| \sqrt{n-1} R(n) B_n^T \right\|^2 = \left| \sqrt{n-1} R^{(1)}(n) \right|^2 \left( \sqrt{n} \ell_\xi(n) \sqrt{\operatorname{Var}(b_1^{(1)} \delta + b_1^{(2)} \varepsilon)} \right)^{-2}$$
$$+ \left\| \sqrt{n-1} R^{(2, d)}(n) n^{-1/2} \left( \operatorname{Cov} J_0^{(2, d)} \right)^{-T/2} \right\|^2.$$

This completes the proof of (3.48).

For establishing (3.49), it suffices to prove that, as $n \to \infty$,

$$\left( (n-1) V_{JJ} \right)^{T/2} B_n^T \xrightarrow{P} I_d, \tag{3.53}$$

and

$$B_n^{-T} \left( (n-1) V_{K(n)K(n)} \right)^{-T/2} \xrightarrow{P} I_d, \tag{3.54}$$

with matrix $B_n$ given by (3.50). Since

$$\left( (n-1) V_{JJ} \right)^{T/2} B_n^T = \left( (n-1) V_{JJ} \right)^{-1/2} B_n^{-1} \left( B_n (n-1) V_{JJ} B_n^T \right),$$

then (3.53) follows from the convergence in (3.51) via (4.5). Similarly, (3.54) is a consequence of (4.5) and

$$B_n (n-1) V_{K(n)K(n)} B_n^T \xrightarrow{P} I_d, \quad n \to \infty. \tag{3.55}$$

The rest of the proof is concerned with verifying (3.55).

Defining vector

$$Q_i(n) := K_i(n) - J_i, \tag{3.56}$$

we have

$$V_{K(n)K(n)} = V_{JJ} + V_{JQ(n)} + V_{Q(n)J} + V_{Q(n)Q(n)}. \tag{3.57}$$

Due to (49) of Lemma 8 in [17] (condition (24) in [17] with $b_j$ as in (3.45) in place of $e$ of (37) in [17] is satisfied on account of (3.25) and (3.29)), as $n \to \infty$,

$$V_{Q^{(j)}(n)Q^{(j)}(n)} = o_P(1), \quad \text{for all } 1 \leq j \leq d, \tag{3.58}$$

where $Q_i^{(j)}(n)$ is the $j^{\text{th}}$ component of $Q_i(n)$. The latter and the Cauchy-Schwarz inequality applied to each off-diagonal entry of matrix $V_{Q(n)Q(n)}$ yeild that $V_{Q(n)Q(n)} \xrightarrow{P} O$ (zero matrix), and therefore,

$$B_n(n-1) V_{Q(n)Q(n)} B_n^T \xrightarrow{P} O, \tag{3.59}$$

with matrix $B_n$ of (3.50), as $n \to \infty$. Now, in view of (3.51), (3.57) and (3.59), convergence in (3.55) is a consequence of

$$B_n(n-1) \left( V_{JQ(n)} + V_{Q(n)J} \right) B_n^T \xrightarrow{P} O, \quad n \to \infty, \tag{3.60}$$

which, in turn, follows from (3.51), (3.59) and the Cauchy-Schwarz inequality applied to each of the entries of the converging matrix product in (3.60). While the latter implication is easily seen when $\operatorname{Var} \xi < \infty$ and/or $b_1^{(1)} = b_1^{(2)} = 0$, otherwise, namely under $\operatorname{Var} \xi = \infty$ and $|b_1^{(1)}| + |b_1^{(2)}| > 0$, one should first obtain the explicit forms of the converging matrix products in (3.51), (3.59) and (3.60) via the lines in (3.15)–(3.18). □

*Proof of the* (b) *part of Theorem* 3.4. The proof follows the same scheme as the one used for the proof of the (a) part of Theorem 3.4. Instead of (3.47) we now have

$$\begin{cases} R^{(1)}(n) = o_P(1), & \text{if } \lim_{n\to\infty} \overline{\xi^2} < \infty \text{ and/or } b_1^{(1)} = b_1^{(2)} = 0, \\ \dfrac{R^{(1)}(n)}{\sqrt{\sum_{i=1}^n (\xi_i - \overline{\xi})^2/n}} = o_P(1), & \text{if } \lim_{n\to\infty} \overline{\xi^2} = \infty \text{ and } |b_1^{(1)}| + |b_1^{(2)}| > 0, \\ R^{(j)}(n) = o_P(1), & \text{for all } 2 \leq j \leq d, \end{cases} \quad (3.61)$$

which follows from (3.53) in [16] and (3.25). In place of $B_n$ of (3.50), we choose here matrix

$$B_n = \begin{cases} \Sigma^{-1/2} n^{-1/2}, & \text{if } \lim_{n\to\infty} \overline{\xi^2} < \infty \text{ and/or } b_1^{(1)} = b_1^{(2)} = 0, \\ \Sigma^{-1/2} D_n, & \text{otherwise}, \end{cases} \quad (3.62)$$

where matrix $\Sigma > 0$ is as in (3.2) that reads for $Z_i(n)$ of (3.31), while matrix $D_n$ is defined in (3.37). Then, convergence in (3.51) with $B_n$ of (3.62) is due to having (3.31), (3.4) and (3.40). As to the validity of (3.52) in the context of the FEIVM (1.1), with $B_n$ of (3.62), it is based on (3.61) and convergence

$$\frac{\sum_{i=1}^n (\xi_i - \overline{\xi})^2}{\sum_{i=1}^n (\xi_i - m)^2} \xrightarrow{P} 1, \quad n \to \infty, \quad (3.63)$$

which amounts to (3.28) in [16]. Thus, (3.48) in the FEIVM is verified via the appropriate versions of (3.51) and (3.52).

Now, we are to argue (3.49) in the present context. In fact, we only need to show (3.59) with $B_n$ of (3.62). By (3.54) in [16] and (3.25), for $Q_i(n)$ as in (3.56),

$$\begin{cases} V_{Q^{(1)}(n)Q^{(1)}(n)} = o_P(1), & \text{if } \lim_{n\to\infty} \overline{\xi^2} < \infty \text{ and/or } b_1^{(1)} = b_1^{(2)} = 0, \\ \dfrac{nV_{Q^{(1)}(n)Q^{(1)}(n)}}{\sum_{i=1}^n (\xi_i - \overline{\xi})^2} = o_P(1), & \text{if } \lim_{n\to\infty} \overline{\xi^2} = \infty \text{ and } |b_1^{(1)}| + |b_1^{(2)}| > 0, \\ V_{Q^{(j)}(n)Q^{(j)}(n)} = o_P(1), & \text{for all } 2 \leq j \leq d. \end{cases} \quad (3.64)$$

Combining (3.63) and (3.64), we conclude (3.59) with $B_n$ of (3.62). □

**Remark 3.4.** The use of Theorem 3.4 can also go beyond the studies of the present paper, since Theorem 3.4 is suitable for establishing joint CLT's for estimators, other than those in (1.3)–(1.5), that are also appropriately based on vector $(\overline{y}, \overline{x}, S_{yy}, S_{xy}, S_{xx})$ in the SEIVM (1.1) and FEIVM (1.1).

**Remark 3.5.** We note that, via [14, 15, 16, 17], Theorem 3.3 can be adapted for the no-intercept versions of the SEIVM (1.1) and FEIVM (1.1), where the intercept $\alpha$ is known to be zero, and respective CLT's for estimators that are appropriately based on vector $(\overline{y}, \overline{x}, \overline{y^2}, \overline{xy}, \overline{x^2})$ can thus be established.

### 3.4. Proof of Theorem 2.1

The proof of Theorem 2.1 is given for the MLSE's $\widehat{\beta}_{3n}$ and $\widehat{\alpha}_{3n}$, and the MME $\widehat{\gamma}_{3n}$ only, simultaneously for the SEIVM (1.1) and FEIVM (1.1). The corresponding CLT's for $(\widehat{\beta}_{jn}, \widehat{\alpha}_{jn}, \widehat{\gamma}_{jn})$ for $j = 1$ and $2$ can be established in similar ways and thus, the respective proofs are omitted here.

*Proof of the* (a) *part of Theorem* 2.1. From the proofs of Theorems 1–3 in [17], and the proofs of the (a) parts of Theorems 2.2 and 2.3 in [16], and also by following the lines of the proofs of Theorems 1.1.2c and 2.1.2c in [15], we have

$$
\begin{aligned}
&\sqrt{n}\, U(3,n)(\widehat{\beta}_{3n} - \beta) = \sqrt{n}\,\overline{u(3,n,\beta)}, \\
&\sqrt{n}\,(\widehat{\alpha}_{3n} - \alpha) = \sqrt{n}\,\overline{v'(3,n,\beta)} + o_P(1), \\
&\sqrt{n}\, L(3,n)(\widehat{\gamma}_{3n} - \gamma) = \sqrt{n}\,\overline{w(3,n,\beta)} + o_P(1),
\end{aligned}
\tag{3.65}
$$

where

$$
v'_i(3,n,\beta) = 
\begin{cases}
y_i - \alpha - \beta x_i - \dfrac{m}{M-m^2}\, u_i(3,n,\beta), & \text{if } E\,\xi^2 = M < \infty \text{ (in the SEIVM)}, \\
& \quad \text{or } \lim_{n\to\infty} \overline{\xi^2} = M < \infty \text{ (in the FEIVM)}, \\
y_i - \alpha - \beta x_i & , \text{ if } \mathrm{Var}\,\xi = \infty, \text{ or } \lim_{n\to\infty} \overline{\xi^2} = \infty,
\end{cases}
\tag{3.66}
$$

with $m = E\,\xi$ (in the SEIVM), or $m = \lim_{n\to\infty} \overline{\xi}$ (in the FEIVM). Put

$$
K_i(n) = \Big(u_i(3,n,\beta), v'_i(3,n,\beta), w(3,n,\beta)\Big). \tag{3.67}
$$

Then from (3.65),

$$
\begin{aligned}
&\sqrt{n}\Big(U(3,n)(\widehat{\beta}_{3n}-\beta),\, \widehat{\alpha}_{3n}-\alpha,\, L(3,n)(\widehat{\gamma}_{3n}-\gamma)\Big) V_{z(3,n,\beta)z(3,n,\beta)}^{-T/2} \\
&=: \Big(\sqrt{n}\,\overline{K(n)} + \rho(n)\Big) V_{z(3,n,\beta)z(3,n,\beta)}^{-T/2} \\
&= \Big(St_n(K(n)) + \rho(n) V_{K(n)K(n)}^{-T/2}\Big) V_{K(n)K(n)}^{T/2} V_{z(3,n,\beta)z(3,n,\beta)}^{-T/2},
\end{aligned}
\tag{3.68}
$$

with

$$
\|\rho(n)\| = o_P(1). \tag{3.69}
$$

Next, we check the conditions of Theorem 3.4 for having $St_n(K(n)) \xrightarrow{\mathcal{D}} N(0, I_d)$, $n \to \infty$. For $K_i(n)$ as in (3.67), the corresponding vectors $c_1$, $c_2$ and $c_3$ as in (3.43) are

$$
\begin{aligned}
c_1 &= (0, 0, 0, 1, -\beta), \\
c_2 &= \begin{cases} (1, -\beta, 0, 0, 0) - \dfrac{m}{M - m^2} c_1, & \text{if } E\,\xi^2 = M < \infty, \\ & \text{or } \lim_{n\to\infty} \overline{\xi^2} = M < \infty, \\ (1, -\beta, 0, 0, 0), & \text{if } \operatorname{Var}\xi = \infty, \text{ or } \lim_{n\to\infty} \overline{\xi^2} = \infty, \end{cases} \\
c_3 &= (0, 0, 1, -2\beta, \beta^2),
\end{aligned} \quad (3.70)
$$

while vectors $b_1$, $b_2$ and $b_3$ as in (3.45) are

$$
\begin{aligned}
b_1 &= (1, -\beta, 0, 0, 1, 0, -\beta), \\
b_2 &= \begin{cases} \left(-\dfrac{m}{M-m^2}, \dfrac{m\beta}{M-m^2}, 1, -\beta, -\dfrac{m}{M-m^2}, 0, \dfrac{m\beta}{M-m^2}\right), \\ \qquad \text{if } E\,\xi^2 = M < \infty, \text{ or } \lim_{n\to\infty} \overline{\xi^2} = M < \infty, \\ (0, 0, 1, -\beta, 0, 0, 0), \text{ if } \operatorname{Var}\xi = \infty, \text{ or } \lim_{n\to\infty} \overline{\xi^2} = \infty, \end{cases} \\
b_3 &= (0, 0, 0, 0, -2\beta, 1, \beta^2).
\end{aligned} \quad (3.71)
$$

It is not hard to see that $c_1$, $c_2$, $c_3$ of (3.70) satisfy (3.44), while $b_1$, $b_2$, $b_3$ of (3.71) obey (3.25). Hence, we only need to check (3.29) and (3.30) with the vectors $b_1$, $b_2$, $b_3$ of (3.71). In fact, according to Remark 3.3, it suffices to verify (3.29) only. Since $|b_1^{(1)}| + |b_1^{(2)}| > 0$, (3.29) reads here as

$$
\begin{cases} \text{vector } (\langle \zeta_0, b_1 \rangle, \langle \zeta_0, b_2 \rangle, \langle \zeta_0, b_3 \rangle) \text{ is full}, & \text{if } \operatorname{Var}\xi < \infty, \\ \text{vector } (\langle \zeta_0, b_2 \rangle, \langle \zeta_0, b_3 \rangle) \text{ is full}, & \text{if } \operatorname{Var}\xi = \infty. \end{cases} \quad (3.72)
$$

Suppose first that $\operatorname{Var}\xi = \infty$. If vector $(\langle \zeta_0, b_2 \rangle, \langle \zeta_0, b_3 \rangle)$ is not full, then there exists vector $u = (u_2, u_3) \in \mathbb{R}^2$, $\|u\| = 1$, such that

$$ u_2 \langle \zeta_0, b_2 \rangle + u_3 \langle \zeta_0, b_3 \rangle = 0 \quad \text{almost surely (a.s.).} \quad (3.73) $$

If $u_3 = 0$, then (3.73) reduces to $\langle \zeta_0, b_2 \rangle = 0$ a.s., or to $\operatorname{Var}\langle \zeta_0, b_2 \rangle = \langle (1, -\beta)\Gamma, (1, -\beta) \rangle = 0$, and the latter equality contradicts positive definiteness of $\Gamma$ as in **(A)**. Assume now that $u_3 \neq 0$ in (3.73). Then it can be shown that (3.73) is equivalent to $(u_2/(2u_3) + (\delta - \beta\varepsilon))^2 = u_2^2/(4u_3^2) + E(\delta - \beta\varepsilon)^2$ a.s., or to $\delta - \beta\varepsilon = \pm\sqrt{u_2^2/(4u_3^2) + E(\delta - \beta\varepsilon)^2} - u_2/(2u_3)$ a.s. that means that $\delta - \beta\varepsilon$ would have to be discretely distributed. However, this would violate **(B)**.

Suppose now that $\operatorname{Var}\xi < \infty$. We are to prove that the vector $(\langle \zeta_0, b_1 \rangle, \langle \zeta_0, b_2 \rangle, \langle \zeta_0, b_3 \rangle)$ is full. Consider r.v.

$$ \psi = \langle \zeta_0, u_1 b_1 + u_2 b_2 + u_3 b_3 \rangle \quad (3.74) $$

for all real $u_1$, $u_2$, $u_3$, where $\|(u_1, u_2, u_3)\| = 1$. We seperate the following four cases in terms of $u_1$ and $u_3$.

1) If $u_1 = u_3 = 0$, then $\operatorname{Var}\psi = \operatorname{Var}\langle \zeta_0, u_2 b_2 \rangle = u_2^2 \operatorname{Var}\langle \zeta_0, b_2 \rangle$. If $m \neq 0$, then $|b_2^{(1)}| + |b_2^{(2)}| > 0$ and $\operatorname{Var}\psi > 0$ by Remark 3.3. Otherwise, when $m = 0$, $\operatorname{Var}\psi = \langle (1, -\beta)\Gamma, (1, -\beta) \rangle > 0$ since $\Gamma > 0$ by **(A)**.

2) Assume that $u_1 = 0$, but $u_3 \neq 0$. Provided that $u_2 \neq 0$ and $m \neq 0$, vector $(u_2 b_2 + u_3 b_3)^{(1,2)} = u_2 b_2^{(1,2)} = -u_2 m(M-m^2)^{-1} b_1^{(1,2)}$ is nonzero and hence, for $\psi$ of (3.74), $\text{Var}\,\psi > 0$ by Remark 3.3. If $u_2 = 0$, then $\text{Var}\,\psi = u_3^2 \text{Var}((\delta - \beta\varepsilon)^2 - E(\delta-\beta\varepsilon)^2) > 0$, because the situation $\delta - \beta\varepsilon = \pm\sqrt{E(\delta-\beta\varepsilon)^2}$ a.s. is prevented via **(B)** ($\delta - \beta\varepsilon$ cannot be discretely distributed). If $m = 0$, then the expression for $b_2$ is as in the case $\text{Var}\,\xi = \infty$ (cf. (3.71)), and the proof of that $\text{Var}\,\psi > 0$ goes exactly like that of disproving (3.73) under $\text{Var}\,\xi = \infty$.

3) Let $u_1 \neq 0$ and $u_3 = 0$. If $u_1 b_1^{(1,2)} + u_2 b_2^{(1,2)} = (u_1 - u_2 m(M-m^2)^{-1}) b_1^{(1,2)}$ is not a zero vector, then from Remark 3.3, $\text{Var}\,\psi = \text{Var}\langle \zeta_0, u_1 b_1 + u_2 b_2 \rangle > 0$, while otherwise, namely, equivalently if $u_1 = u_2 m(M-m^2)^{-1}$, then $u_1 b_1^{(5,7)} + u_2 b_2^{(5,7)}$ is a zero vector and, due to positive definiteness of $\Gamma$ assumed in **(A)**,

$$\text{Var}\,\psi = \text{Var}\langle \zeta_0, u_1 b_1 + u_2 b_2 \rangle = \text{Var}\langle \zeta_0^{(3,7)}, (u_1 b_1 + u_2 b_2)^{(3,7)} \rangle$$
$$= \text{Var}\langle \zeta_0^{(3,4)}, u_2 b_2^{(3,4)} \rangle = u_2^2 \langle (1, -\beta)\Gamma, (1, -\beta)\rangle > 0.$$

4) Finally, consider the case of $u_1 \neq 0$ and $u_3 \neq 0$. If vector $u_1 b_1^{(1,2)} + u_2 b_2^{(1,2)} + u_3 b_3^{(1,2)}$ is nonzero, then $\text{Var}\,\psi > 0$ by Remark 3.3. Otherwise, namely, if $u_1 b_1^{(1,2)} + u_2 b_2^{(1,2)} + u_3 b_3^{(1,2)} = u_1 b_1^{(1,2)} + u_2 b_2^{(1,2)} = (u_1 - u_2 m(M-m^2)^{-1}) b_1^{(1,2)}$ is a zero vector, then so is vector $u_1 b_1^{(5,7)} + u_2 b_2^{(5,7)} = (u_1 - u_2 m(M-m^2)^{-1}) b_1^{(5,7)}$ and $\text{Var}\,\psi = \text{Var}\langle \zeta_0^{(3,7)}, (u_1 b_1 + u_2 b_2 + u_3 b_3)^{(3,7)} \rangle = \text{Var}(u_2 \langle \zeta_0^{(3,4)}, b_2^{(3,4)}\rangle + u_3 \langle \zeta_0^{(5,7)}, b_3^{(5,7)} \rangle)$. Note that $(\langle \zeta_0^{(3,4)}, b_2^{(3,4)} \rangle, \langle \zeta_0^{(5,7)}, b_3^{(5,7)} \rangle) = (\langle \zeta_0, \tilde{b}_2 \rangle, \langle \zeta_0, b_3 \rangle)$, with vector $\tilde{b}_2$ equal to $b_2$ as in case of $\text{Var}\,\xi = \infty$ (cf. (3.71)). Also, $|u_2| + |u_3| > 0$. Consequently, using the just obtained expression for $\text{Var}\,\psi$, we conclude that, on account of that (3.73) fails, $\text{Var}\,\psi > 0$.

This completes verification of (3.72).

Consequently, by Theorem 3.4, $St_n(K(n)) \xrightarrow{\mathcal{D}} N(0, I_d)$, as $n \to \infty$. Hence, in view of (3.68), the proof reduces to showing that, as $n \to \infty$,

$$\|\rho(n) V_{K(n)K(n)}^{-T/2}\| = o_P(1) \tag{3.75}$$

and

$$V_{K(n)K(n)}^{T/2} V_{z(3,n,\beta)z(3,n,\beta)}^{-T/2} \xrightarrow{P} I_3. \tag{3.76}$$

Similarly to the proof of (3.48), (3.75) is on account of (3.55) (proved both for the SEIVM and FEIVM versions of (1.1)) and

$$\|\sqrt{n-1}\rho(n) B_n^T\| = o_P(1), \tag{3.77}$$

where $\rho(n)$ is as in (3.69), matrix $B_n$ is as in (3.50) in the SEIVM, with vector $J_0 = (\langle \zeta_0, b_1 \rangle, \langle \zeta_0, b_2 \rangle, \langle \zeta_0, b_3 \rangle)$, or as in (3.62) in the FEIVM, with matrix $\Sigma$ of (3.2) corresponding to vectors $Z_i(n)$ of (3.31) that are based on $J_i = (\langle \zeta_i, b_1 \rangle, \langle \zeta_i, b_2 \rangle, \langle \zeta_i, b_3 \rangle)$, where the vectors $b_1$, $b_2$ and $b_3$ are as in (3.71).

As to the convergence in (3.76), it suffices to show that, as $n \to \infty$,

$$\left((n-1)V_{K(n)K(n)}\right)^{T/2} B_n^T \xrightarrow{P} I_3 \tag{3.78}$$

and

$$B_n^{-T} \left((n-1)V_{z(3,n,\beta)z(3,n,\beta)}\right)^{-T/2} \xrightarrow{P} I_3, \tag{3.79}$$

with matrix $B_n$ as above. Since

$$\left((n-1)V_{K(n)K(n)}\right)^{T/2} B_n^T$$
$$= \left((n-1)V_{K(n)K(n)}\right)^{-1/2} B_n^{-1} \left(B_n(n-1)V_{K(n)K(n)}B_n^T\right),$$

(3.78) is a consequence of (3.55) (proved both for the SEIVM and FEIVM) and (4.5). In the same manner, by (4.5), (3.79) is due to

$$B_n(n-1)V_{z(3,n,\beta)z(3,n,\beta)}B_n^T \xrightarrow{P} I_3, \quad n \to \infty. \tag{3.80}$$

By writing

$$V_{z(3,n,\beta)z(3,n,\beta)} = V_{K(n)K(n)} + V_{K(n)\Delta(n)} + V_{\Delta(n)K(n)} + V_{\Delta(n)\Delta(n)},$$

where vectors

$$\Delta_i(n) = z_i(3,n,\beta) - K_i(n) = (0, v_i(3,n,\beta) - v'_i(3,n,\beta), 0), \tag{3.81}$$

with $v_i(3,n,\beta)$ of (2.3) and $v'_i(3,n,\beta)$ of (3.66), and by using similar arguments to those for (3.55), (3.80) results from (3.55) and

$$B_n(n-1)V_{\Delta(n)\Delta(n)}B_n^T \xrightarrow{P} 0, \quad n \to \infty. \tag{3.82}$$

As to the validity of (3.82), the lines in [17] that are right below (64) and (67) imply that in the SEIVM $V_{\Delta^{(2)}(n)\Delta^{(2)}(n)} = o_P(1)$, $n \to \infty$. Mutatis mutandis, the latter holds true also in the FEIVM. This leads to (3.82). □

*Proof of the* (b) *part of Theorem 2.1.* If, as $n \to \infty$,

$$V_{z(3,n,\beta)z(3,n,\beta)}^{T/2} V_{z(3,n,\widehat{\beta}_{jn})z(3,n,\widehat{\beta}_{jn})}^{-T/2} \xrightarrow{P} I_3, \tag{3.83}$$

then the (b) part of Theorem 2.1 follows from its (a) part. As $n \to \infty$,

$$\left((n-1)V_{z(3,n,\beta)z(3,n,\beta)}\right)^{T/2} B_n^T \xrightarrow{P} I_3 \tag{3.84}$$

and

$$B_n^{-T} \left((n-1)V_{z(3,n,\widehat{\beta}_{jn})z(3,n,\widehat{\beta}_{jn})}\right)^{-T/2} \xrightarrow{P} I_3 \tag{3.85}$$

guarantee (3.83), where matrix $B_n$ is given by the appropriate versions of (3.50) and (3.62) in the SEIVM and FEIVM respectively, just like in the proof of the (a) part of Theorem 2.1.

Using (4.5), (3.80) and representation

$$((n-1)V_{z(3,n,\beta)z(3,n,\beta)})^{T/2} B_n^T$$
$$= ((n-1)V_{z(3,n,\beta)z(3,n,\beta)})^{-1/2} B_n^{-1} \left(B_n(n-1)V_{z(3,n,\beta)z(3,n,\beta)}B_n^T\right),$$

we conclude (3.84).

Since

$$V_{z(3,n,\widehat{\beta}_{jn})z(3,n,\widehat{\beta}_{jn})} = V_{z(3,n,\beta)z(3,n,\beta)} + V_{z(3,n,\beta)p(n)} + V_{p(n)z(3,n,\beta)} + V_{p(n)p(n)},$$

with vectors

$$p_i(n) = z(3,n,\widehat{\beta}_{jn}) - z_i(3,n,\beta), \tag{3.86}$$

then by similar arguments to those for (3.55), convergence in (3.85) is a consequence of (3.80) and

$$B_n(n-1)V_{p(n)p(n)}B_n^T \xrightarrow{P} 0, \quad n \to \infty. \tag{3.87}$$

Now we argue the validity of (3.87). From the lines of the proofs of Theorems 2 and 3 in [17], the (b) parts of Theorems 2.2 and 2.3 in [16] and Theorems 1.1.3c and 2.1.3c in [15], we have

$$\begin{cases} \dfrac{V_{p^{(1)}(n)p^{(1)}(n)}}{\ell_\xi^2(n)} = o_P(1), & \text{in the SEIVM,} \\ \dfrac{nV_{p^{(1)}(n)p^{(1)}(n)}}{\sum_{i=1}^n (\xi_i - \overline{\xi})^2} = o_P(1), & \text{in the FEIVM,} \\ V_{p^{(2)}(n)p^{(2)}(n)} = o_P(1), V_{p^{(3)}(n)p^{(3)}(n)} = o_P(1) \end{cases} \tag{3.88}$$

and hence, also (3.87). □

*Proof related to Remark* 2.3. In view of (3.83), it suffices to show that under (2.6) (in the SEIVM) and (2.7) (in the FEIVM), subject to the additional normality model conditions in [3] and comments on **(B)** in Remarks 2.8 and 2.9, the (a) part of Theorem 2.1 implies corresponding $\sqrt{n}$−asymptotic normality of $(\widehat{\beta}_{jn}, \widehat{\alpha}_{jn}, \widehat{\gamma}_{jn})$, $1 \le j \le 3$, as in [3, 7, 8]. Due to similar arguments, we consider here the case of $j = 3$ only.

From (3.84) with $B_n = n^{-1/2} \Big( \text{Cov}(\langle \zeta_0, b_1 \rangle, \langle \zeta_0, b_2 \rangle, \langle \zeta_0, b_3 \rangle) \Big)^{-1/2}$ in the SEIVM and $B_n = \Big( \sum_{i=1}^n \text{Cov}(\langle \zeta_i, b_1 \rangle, \langle \zeta_i, b_2 \rangle, \langle \zeta_i, b_3 \rangle) \Big)^{-1/2}$ in the FEIVM, where vectors $b_1$, $b_2$, $b_3$ are as in (3.71), and from the corresponding (3.29) and (3.30) that have been checked in the (a) part of Theorem 2.1 via assuming **(B)**,

$$V_{z(3,n,\beta)z(3,n,\beta)}^{T/2} \xrightarrow{P} A^{T/2} :=$$
$$\begin{cases} \Big( \text{Cov}(\langle \zeta_0, b_1 \rangle, \langle \zeta_0, b_2 \rangle, \langle \zeta_0, b_3 \rangle) \Big)^{T/2}, & \text{in the SEIVM,} \\ \Big( \lim_{n \to \infty} n^{-1} \sum_{i=1}^n \text{Cov}(\langle \zeta_i, b_1 \rangle, \langle \zeta_i, b_2 \rangle, \langle \zeta_i, b_3 \rangle) \Big)^{T/2}, & \text{in the FEIVM.} \end{cases}$$
$$\tag{3.89}$$

Next, as $n \to \infty$, since $U(3,n) = S_{xx} - \theta \overset{P}{\to} M - m^2 > 0$, with $M$ and $m$ as in (2.6) or (2.7), then

$$\left(\mathrm{diag}(U(3,n),1,L(3,n))\right)^{-1} \overset{P}{\to} \left(\mathrm{diag}(M-m^2,1,1)\right)^{-1} := D^{-1}, \quad n \to \infty. \tag{3.90}$$

Finally, the (a) part of Theorem 2.1 with $j = 3$ that reads as

$$\sqrt{n}\left(U(3,n)(\widehat{\beta}_{3n} - \beta), \widehat{\alpha}_{3n} - \alpha, L(3,n)(\widehat{\gamma}_{3n} - \gamma)\right)V_{z(3,n,\beta)z(3,n,\beta)}^{-T/2}$$
$$= \sqrt{n}\left(\widehat{\beta}_{3n} - \beta, \widehat{\alpha}_{3n} - \alpha, \widehat{\gamma}_{3n} - \gamma\right)\mathrm{diag}(U(3,n),1,L(3,n))V_{z(3,n,\beta)z(3,n,\beta)}^{-T/2}$$
$$\overset{\mathcal{D}}{\to} N(0, I_3), \quad n \to \infty,$$

and convergence in (3.89) and (3.90) imply that

$$\sqrt{n}\left(\widehat{\beta}_{3n} - \beta, \widehat{\alpha}_{3n} - \alpha, \widehat{\gamma}_{3n} - \gamma\right) \overset{\mathcal{D}}{\to} N(0, D^{-1}AD^{-1}), \quad n \to \infty, \tag{3.91}$$

where matrix $D^{-1}AD^{-1}$ reduces to the covariance matrix of the corresponding asymptotically normal distribution that is obtained in [3], subject to the therein assumed additional normality model conditions spelled out in Remark 2.3 that replace our weaker assumption **(B)** that was used to conclude (3.91) here (for a summary on **(B)** we refer to Remarks 2.8 and 2.9). □

## 4. Appendix: some results on Studentization of random vectors by a matrix and the generalized domain of attraction of the multivariate normal law

This appendix of some well-known results on recent advances in Studentization of random vectors by a matrix and the generalized domain of attraction of the multivariate normal law (GDAN) is provided here for convenience in references that are required for reading Section 3 of this paper.

Motivated by the importance of having matrix Studentized CLT's for random vectors converging in distribution to a spherically symmetric random vector $Z$, which is such that all its Euclidean inner products $\langle Z, u \rangle$ with a deterministic vector $u$ of unit Euclidean norm have the same dirtibution coinciding with that of each single component of $Z$, Vu, Maller and Klass [22] established a rather general result of such a nature. Their result that follows herewith is essentially a general recipe of matrix "Slutskying" for random vectors.

**Theorem 4.1.** *For $n \geq 1$, let $C_n$ be $d \times d$ real invertible nonstochastic matrices, $V_n$ be $d \times d$ real symmetric stochastic matrices and $S_n$ be $1 \times d$ stochastic vectors. If, as $n \to \infty$,*

$$C_n V_n C_n^T \overset{P}{\to} I_d, \tag{4.1}$$

*with unit $d \times d$ matrix $I_d$, and*

$$S_n C_n^T \overset{\mathcal{D}}{\to} Z, \tag{4.2}$$

where $Z$ is a spherically symmetric random vector in $\mathbb{R}^d$, then

$$S_n V_n^{-T/2} \xrightarrow{\mathcal{D}} Z. \tag{4.3}$$

**Remark 4.1.** We note that (4.1) amounts to saying that each entry of matrix $C_n V_n C_n^T$ converges in probability to the corresponding entry of $I_d$. As noted in [22], (4.1) implies that $V_n > 0$ on sets whose probabilities converge to one, as $n \to \infty$. Hence, the Cholesky and the symmetric positive definite square roots of $V_n$ and, consequently, $V_n^{-T/2}$ in (4.3) are well-defined.

**Remark 4.2.** This is a technical remark related to Theorem 4.1 and frequently used in Section 3. It is noted in [22] that if instead of (4.1) and (4.2) one assumes, as $n \to \infty$,

$$A_n^{-1/2} V_n A_n^{-T/2} \xrightarrow{P} I_d \quad \text{and} \quad S_n A_n^{-T/2} \xrightarrow{\mathcal{D}} Z, \tag{4.4}$$

with $d \times d$ real matrices $A_n > 0$, $n \geq 1$, then spherical symmetry of $Z$ is not required for the conclusion in (4.3). Moreover, if (4.4) is assumed to begin with in [22], then the proof of corresponding (4.3) there reduces to showing that, as $n \to \infty$,

$$A_n^{-1/2} V_n A_n^{-T/2} \xrightarrow{P} I_d \text{ implies } A_n^{T/2} V_n^{-T/2} = \left(V_n^{-1/2} A_n^{1/2}\right)^T \xrightarrow{P} I_d. \tag{4.5}$$

A definition of the generalized domain of attraction of the $d$-variate normal law, denoted here by GDAN in view of the previously used DAN as in **(S1)**, reads as follows.

**Definition 4.1.** *Let $\{Z, Z_i, i \geq 1\}$ be i.i.d. random vectors in $\mathbb{R}^d$. $Z$ is said to belong to* GDAN *if there exist nonstochastic sequences of vectors $a_n$ and $d \times d$ matrices $B_n$, such that*

$$\left(\sum_{i=1}^n Z_i - a_n\right) B_n^T \xrightarrow{\mathcal{D}} N(0, I_d), \quad as \ \ n \to \infty. \tag{4.6}$$

**Remark 4.3.** It was shown in Maller [13] that if (4.6) holds, then $E\|Z\|^\alpha < \infty$ for all $0 \leq \alpha < 2$ and $a_n$ can be taken as $nEZ$, while norming matrix $B_n$ is invertible for large enough $n$ and may be chosen to be symmetric. Also, $B_n \to 0$, as $n \to \infty$. As a general fact, it is also known that (4.6) implies that $Z$ is full (cf. Lemma 3.3.3 in [18]).

Below we review some known equivalent characterizations of GDAN that are found in Theorem 1.1 of Maller [13] and are most relevant to the aims of this paper.

**Theorem 4.2.** *Let $\{Z, Z_i, i \geq 1\}$ be i.i.d. random vectors in $\mathbb{R}^d$ having a full distribution. The following statements are equivalent:*

(a) $Z \in$ GDAN;

(b) *there exist nonstochastic square matrices $B_n$ such that*
$$B_n \sum_{i=1}^{n}(Z_i - \overline{Z})^T(Z_i - \overline{Z})B_n^T \xrightarrow{P} I_d, \quad \text{as } n \to \infty;$$

(c) $\sup\limits_{u \in \mathbb{R}^d,\ \|u\|=1} \dfrac{x^2 P\Big(|\langle Z, u\rangle| > x\Big)}{E\Big(\langle Z, u\rangle^2 \mathbb{1}_{\{|\langle Z,u\rangle|\leq x\}}\Big)} \to 0, \quad \text{as } x \to \infty.$

**Remark 4.4.** Equivalence of (a) and (c) in Theorem 4.2 under $EZ = 0$ is one of the main results of the seminal paper by Hahn and Klass [10] that has stimulated intensive studies of GDAN. The (b) part of Theorem 4.2 can be viewed as an analogue of the following result for DAN that was rediscovered by Maller [12] to conclude a CLT for self-normalized partial sums of i.i.d.r.v.'s from DAN. It is essentially a variation of Theorems 4 and 5 on pp. 143–144 in Gnedenko and Kolmogorov [9]. Accordingly, if $\{Z, Z_i, i \geq 1\}$ are i.i.d.r.v.'s in DAN, then
$$\sum_{i=1}^{n}(Z_i - E\,Z)^2/b_n^2 \xrightarrow{P} 1, \quad n \to \infty, \tag{4.7}$$

where $b_n$ is such that $\sum_{i=1}^{n}(Z_i - E\,Z)/b_n \xrightarrow{D} N(0,1)$, as $n \to \infty$.

**Remark 4.5.** It is easy to see that if $Z \in$ GDAN, then each component $Z^{(j)}$ of $Z$ is in DAN, $1 \leq j \leq d$. Indeed, the (c) part of Theorem 4.2 implies that
$$\frac{x^2 P\Big(|Z^{(j)}| > x\Big)}{E\Big((Z^{(j)})^2 \mathbb{1}_{\{|Z^{(j)}|\leq x\}}\Big)} \to 0, \quad x \to \infty. \tag{4.8}$$

The latter convergence is known as Lévy's necessary and sufficient condition for $Z^{(j)}$ to be in DAN (cf. [11]). On the other hand, assuming that all $Z^{(j)} \in$ DAN is not sufficient alone to guarantee that $Z \in$ GDAN. Indeed, suppose that $EZ = 0$ and all $Z^{(j)}$ are identically distributed and belong to DAN. Then, from (4.7) with $b_n$ chosen to be the same for all sequences $\{Z_i^{(j)}, i \geq 1\}$, $1 \leq j \leq d$,
$$\frac{\sum_{i=1}^{n}\|Z_i\|^2}{d\,b_n^2} \xrightarrow{P} 1, \quad n \to \infty. \tag{4.9}$$

Further, from Remark (ii) on p.193 of [13], via p.236 of [5], (4.9) is equivalent to (4.8) with $\|Z\|$ in place of $|Z^{(j)}|$, and such a form of (4.8) does not imply (c) of Theorem 4.2. However, the class of spherically symmetric random vectors $Z$ is pointed out in Remark (ii) on p.217 of [13] as an exception in this regard, on account of each projection of $Z$ having the same distribution coinciding with that of each single component $Z^{(j)}$ of $Z$, $1 \leq j \leq d$.

For a multivariate Student statistic based on i.i.d. random vectors $\{Z, Z_i, i \geq 1\}$ in GDAN, by Maller [13] and Vu, Maller and Klass [22], the following CLT holds true.

**Theorem 4.3.** *Let $\{Z, Z_i, i \geq 1\}$ be i.i.d. random vectors in $\mathbb{R}^d$ and $Z \in$ GDAN. Then, for the multivariate Student statistic $St_n(Z) = \sqrt{n}\,\overline{Z}\left((n-1)^{-1} \sum_{i=1}^{n}(Z_i - \overline{Z})^T(Z_i - \overline{Z})\right)^{-T/2}$, as $n \to \infty$, $St_n(Z - E\,Z) \xrightarrow{\mathcal{D}} N(0, I_d)$.*

**Remark 4.6.** As noted in Remarks of [22], the proof of Theorem 4.3 goes by Theorem 4.1, via combining (4.6) and the (b) part of Theorem 4.2 that presents a special case of (4.1). It is pointed out in [22] that, though Studentization in Theorems 4.1 and 4.3 can be performed by both the Cholesky and symmetric positive definite square roots, one's preference would depend on the problem in hand. It is also conjectured in [22] that for the purpose of transforming $\sum_{i=1}^{n}(Z_i - EZ)$ so as to have approximately a spherically symmetric distribution, the symmetric positive definite square root is likely a better choice to accomplish this in small samples.

## Acknowledgements

The author is thankful to her Ph.D. thesis supervisor Miklós Csörgő for his advice and comments that helped in preparing the present exposition.